\documentclass[a4paper,12pt,reqno]{amsart}
\usepackage[colorlinks=true]{hyperref}
\usepackage{amssymb}
\usepackage{pstricks}
\usepackage{graphicx}
\usepackage{mathdots}
\allowdisplaybreaks[1]
\newtheorem{theorem}{Theorem}
\newtheorem*{theorem*}{Theorem}
\newcommand{\ASM}{\mathrm{ASM}}
\newcommand{\DPP}{\mathrm{DPP}}
\newcommand{\SVDWBC}{\mathrm{6VDW}}
\newcommand{\SV}{\mathrm{6V}}
\newcommand{\NILP}{\mathrm{NILP}}

\author[R.~E.~Behrend]{Roger E.~Behrend}
\address{R.~E.~Behrend, School of Mathematics, Cardiff University, Cardiff, CF24 4AG, UK}
\email{behrendr@cardiff.ac.uk}
\gdef\s{s}
\title[A doubly-refined enumeration of ASM{\s} and DPP{\s}]{A doubly-refined enumeration of
alternating sign matrices and descending plane partitions}
\author[P.~Di Francesco]{Philippe Di Francesco}
\address{P.~Di Francesco,
Institut de Physique Th\'eorique de Saclay,
CEA/DSM/SPhT, CNRS URA 2306,
C.E.A.-Saclay, F-91191 Gif sur Yvette Cedex, France}
\email{philippe.di-francesco@cea.fr}
\author[P.~Zinn-Justin]{Paul Zinn-Justin}
\address{P.~Zinn-Justin, UPMC Univ.~Paris 6, CNRS UMR 7589, LPTHE,
75252 Paris Cedex, France}
\email{pzinn@lpthe.jussieu.fr}
\thanks{PDF and PZJ acknowledge partial support
from ANR program ``GRANMA'' BLAN08-1-13695.
PZJ is supported in part by ERC grant 278124 ``LIC''}
\keywords{Alternating sign matrices, descending plane partitions, six-vertex model with domain-wall boundary
conditions, nonintersecting lattice paths, Desnanot--Jacobi identity}
\begin{document}
\begin{abstract}
It was shown recently by the authors that, for any $n$, there is equality between
the distributions of certain triplets of statistics on $n\times n$ alternating sign matrices (ASMs)
and descending plane partitions (DPPs) with each part at most~$n$.
The statistics for an ASM~$A$ are the number of generalized inversions in~$A$, the number of~$-1$'s in~$A$ and
the number of~$0$'s to the left of the~$1$ in the first row of~$A$, and the respective statistics for a DPP~$D$ are the number of
nonspecial parts in~$D$, the number of special parts in~$D$ and the number of~$n$'s in~$D$.
Here, the result is generalized to include a fourth statistic for each type of object,
where this is the number of~$0$'s to the right of the~$1$ in the last row of an ASM,
and the number of $(n-1)$'s plus the number of rows of length~$n-1$ in a DPP.
This generalization is proved using the known equality of the three-statistic generating functions, together with
relations which express each four-statistic generating function in terms of
its three-statistic counterpart.  These relations are obtained by applying the Desnanot--Jacobi
identity to determinantal expressions for the generating functions,
where the determinants arise from standard methods involving
the six-vertex model with domain-wall boundary conditions
for ASMs, and nonintersecting lattice paths for DPPs.
\end{abstract}
\maketitle

\section{Introduction}\label{intro}
Alternating sign matrices (ASMs) and descending plane partitions (DPPs) are combinatorial objects which share some basic enumerative
properties.  (For an introduction to ASMs, DPPs and related subjects see, for example,
the reviews of Bressoud~\cite{Bre99}, Bressoud and Propp~\cite{BrePro99},
Propp~\cite{Pro01}, Robbins~\cite{Rob91} or Zeilberger~\cite{Zei05}.  For a more up-to-date, but briefer,
introduction see~\cite[Sec.~1]{BehDifZin12a}.)

These objects initially appeared in different contexts, but at approximately the same time,
and it was shortly thereafter that the enumerative connection between them was first observed.
In particular, a conjecture was formulated by Mills, Robbins and Rumsey~\cite[Conj.~3]{MilRobRum83},
stating that the distributions of certain triplets of statistics
on $n\times n$ ASMs and on DPPs with each part at most~$n$ are equal for any positive integer~$n$.
Interestingly, at this time, an aspect of the unproved conjecture was helpful to Mills, Robbins and
Rumsey in proving a different conjecture, namely the Macdonald conjecture
for the generating function of cyclically symmetric plane partitions~\cite{MilRobRum82}.

Some special cases of the conjecture relating ASMs and DPPs were subsequently confirmed (see~\cite[Sec.~1.3]{BehDifZin12a} for a
detailed overview of these),
and the full conjecture was proved recently by the present authors~\cite[Thm.~1]{BehDifZin12a}.
Among the previously-proved special cases, probably the best-known is simply that the total numbers of $n\times n$
ASMs and DPPs with each part at most~$n$ are equal for any~$n$.  This follows from formulae
of Andrews~\cite{And79} for DPPs, and of Zeilberger~\cite{Zei96a} or Kuperberg~\cite{Kup96} for ASMs,
which give the numbers of each type of object as
$\prod_{i=0}^{n-1}\frac{(3i+1)!}{(n+i)!}$.

The primary aim of this paper is to generalize the result of~\cite[Thm.~1]{BehDifZin12a},
involving equally-distributed triplets of statistics on ASMs and DPPs,
to a result involving equally-distributed quadruplets of statistics.
The three original statistics for an ASM~$A$ are the number of generalized inversions in~$A$, the number of~$-1$'s in~$A$ and
the number of~$0$'s to the left of the~$1$ in the first row of~$A$, and the respective statistics for a DPP~$D$ with each part at most~$n$
are the number of nonspecial parts in~$D$, the number of special parts in~$D$ and the number of~$n$'s in~$D$.
The additional statistics considered here are the
number of~$0$'s to the right of the~$1$ in the last row of an ASM,
and the number of $(n-1)$'s plus the number of rows of length~$n-1$ in a DPP with each part at most~$n$.
In each case, the first and second statistics depend
on the bulk structure of the object, while the third and fourth statistics depend only on the structure at or near its boundary.

The joint distribution of the fourth ASM statistic with some or all of the other three ASM statistics has
previously been considered.
See Section~\ref{doubrefsec} for an overview of such studies in the literature.
However, it seems that the fourth DPP statistic has not appeared before in the literature, and
that no other previously-studied DPP statistic has the same joint distribution with the
other three DPP statistics as that of this statistic.

In studies of ASMs, DPPs and related objects,
it is common to regard a refined enumeration as having a certain order, where this is based on only
the number of boundary statistics being used.  Hence, the enumeration of this paper is doubly-refined, although it
involves two bulk statistics in addition to its two boundary statistics.

The method of proving the main result of this paper is to introduce four-statistic,
or doubly-refined, generating functions for ASMs and DPPs, and to show that these can be expressed
as identical combinations of the respective three-statistic, or singly-refined, generating functions.
(The expression for the ASM case was previously obtained by Colomo and Pronko~\cite[Eq.~(5.32)]{ColPro05b},~\cite[Eq.~(3.32)]{ColPro06},
but a new proof is given here.)
The required equality of the ASM and DPP doubly-refined generating functions then follows from the equality,
as confirmed in~\cite[Thm.~1]{BehDifZin12a},
of the ASM and DPP singly-refined generating functions.  The identities which give each doubly-refined generating function
in terms of its singly-refined counterpart are obtained by applying a certain form of the Desnanot--Jacobi determinant identity to
determinantal expressions for the generating functions.  The latter expressions are themselves obtained
using standard techniques involving the statistical mechanical six-vertex model with domain-wall boundary conditions and
the Izergin--Korepin formula~\cite{Ize87,Kor82} for the ASM case,
and certain sets of nonintersecting lattice paths and the Lindstr\"{o}m--Gessel--Viennot theorem~\cite{GesVie85,GesVie89,Lin73} for the DPP case.
These techniques also played an essential role in the proof of~\cite[Thm.~1]{BehDifZin12a}.

An outline of the rest of this paper is as follows.  In Section~\ref{defin}, the definitions of ASMs and DPPs,
and of the relevant statistics and associated generating functions, are given.
In Section~\ref{mainres}, the main result involving the equality of the ASM and DPP doubly-refined generating functions,
and the identities expressing each doubly-refined generating function in terms of its singly-refined counterpart, are stated.
In Sections~\ref{ASMpropsect} and~\ref{DPPpropsect}, proofs of the latter identities for the ASM and DPP cases respectively are given, with
the Desnanot--Jacobi determinant identity, which is central to each of these proofs, having been introduced in
Section~\ref{DesJacSec}.  In Section~\ref{discuss}, some further aspects of this work are discussed.

\section{Definitions}\label{defin}
\subsection{ASMs and DPPs}
In this section, the standard definitions of ASMs and DPPs, and of sets of ASMs and DPPs of order $n$, are given.

An ASM, as first defined by Mills, Robbins and Rumsey~\cite{MilRobRum82,MilRobRum83},
is a square matrix in which each entry is~$0$, $1$ or~$-1$, the nonzero entries alternate in sign along each row and column,
and the sum of entries in each row and column is $1$.
It follows that an ASM has a unique 1 in each of its first and last row and column,
and that any permutation matrix is an ASM.

A DPP, as first defined by Andrews~\cite{And79,And80}, is an array of positive integers, called parts,
of the form
\setlength{\unitlength}{5pt}
\begin{equation}\label{DPP}\begin{array}{@{}c@{\:}c@{\:}c@{\:}l@{}}
D_{11}&D_{12}&D_{13}&\begin{picture}(32.1,0)\multiput(0,0)(1,0){32}{.}\end{picture}D_{1,\lambda_1}\\
&D_{22}&D_{23}&\begin{picture}(27.1,0)\multiput(0,0)(1,0){27}{.}\end{picture}D_{2,\lambda_2+1}\\
&&D_{33}&\begin{picture}(22.1,0)\multiput(0,0)(1,0){22}{.}\end{picture}D_{3,\lambda_3+2}\\
&&&\begin{picture}(20,3.3)\multiput(1,0)(-1,1){3}{.}\multiput(19,0)(1,1){3}{.}\end{picture}\\
&&&\hspace*{8pt}D_{tt}\begin{picture}(8.6,0)\multiput(0.5,0)(1,0){8}{.}\end{picture}D_{t,\lambda_t+t-1}\;,\end{array}\end{equation}
where
the parts decrease weakly along rows and
decrease strictly down columns,
and the first parts of the rows $D_{ii}$ and row lengths $\lambda_i$ satisfy
\begin{equation}\label{rowcond}D_{11}>\lambda_1\ge D_{22}>\lambda_2\ge\ldots\ge D_{t-1,t-1}>\lambda_{t-1}\ge D_{tt}>\lambda_t.\end{equation}
The empty array is also regarded as a DPP, and denoted $\emptyset$.

Examples of an ASM and a DPP (as also used in~\cite[Eqs.~(1) \&~(4)]{BehDifZin12a}) are
\begin{equation}\label{ASMDPPEx}
A=\left(\begin{array}{@{}c@{\,\;}c@{\,\;}c@{\,\;}c@{\,\;}c@{\,\;}c@{}}0&0&0&1&0&0\\[-0.5mm]0&1&0&-1&1&0\\[-0.5mm]1&-1&1&0&0&0
\\[-0.5mm]0&0&0&1&0&0\\[-0.5mm]0&1&0&-1&0&1\\[-0.5mm]0&0&0&1&0&0\end{array}\right),\qquad
\begin{array}{@{}c@{\;}c@{\;}c@{\;}c@{\;}c@{\;}c@{}}&6&6&6&5&2\\D=&&4&4&1\\&&&3\end{array}.\end{equation}

For each positive integer $n$, the sets of ASMs and DPPs of order $n$ are defined as
\begin{equation}\label{ASMDPPn}\ASM(n)=\{n\times n\text{ ASMs}\},\qquad
\DPP(n)=\{\text{DPPs with each part}\le n\}.\end{equation}
For example,
\begin{align}
\notag\ASM(3)&=\left\{\!\left(\begin{array}{@{}c@{\;}c@{\;}c@{}}1&0&0\\[-0.7mm]0&1&0\\[-0.7mm]0&0&1\end{array}\right)\!,
\left(\begin{array}{@{}c@{\;}c@{\;}c@{}}0&0&1\\[-0.7mm]0&1&0\\[-0.7mm]1&0&0\end{array}\right)\!,
\left(\begin{array}{@{}c@{\;}c@{\;}c@{}}1&0&0\\[-0.7mm]0&0&1\\[-0.7mm]0&1&0\end{array}\right)\!,
\left(\begin{array}{@{}c@{\;}c@{\;}c@{}}0&0&1\\[-0.7mm]1&0&0\\[-0.7mm]0&1&0\end{array}\right)\!,
\left(\begin{array}{@{}c@{\;}c@{\;}c@{}}0&1&0\\[-0.7mm]1&0&0\\[-0.7mm]0&0&1\end{array}\right)\!,
\left(\begin{array}{@{}c@{\;}c@{\;}c@{}}0&1&0\\[-0.7mm]0&0&1\\[-0.7mm]1&0&0\end{array}\right)\!,
\left(\begin{array}{@{}c@{\:}c@{\;}c@{}}0&1&0\\[-0.7mm]1&-1&1\\[-0.7mm]0&1&0\end{array}\right)\!\right\},\\
\label{ASMDPP3}\DPP(3)&=\left\{\emptyset,\;\begin{array}{@{}c@{\;}c@{}}3&3\\[-0.8mm]&2\end{array},\;2,\;3\;3,\;3,
\;3\;2,\;3\;1\right\}.\end{align}

\subsection{Statistics and generating functions for ASMs and DPPs}\label{stat}
In this section, certain statistics, and associated generating functions, are defined for ASMs and DPPs.
Some simple properties of the generating functions are also identified.

For a given positive integer $n$, statistics for each $A\in\ASM(n)$ are defined as
\begin{align}
\label{nuA}\nu(A)&=\sum_{\substack{1\le i<i'\le n\\1\le j'\le j\le n}}\!A_{ij}\,A_{i'j'},\\[1mm]
\label{muA}\mu(A)&=\text{number of $-1$'s in }A,\\[1mm]
\label{rho1A}\rho_1(A)&=\text{number of 0's to the left of the 1 in the first row of }A,\\[1mm]
\label{rho2A}\rho_2(A)&=\text{number of 0's to the right of the 1 in the last row of }A,\\
\intertext{and statistics for each $D\in\DPP(n)$ are defined as}
\label{nuD}\nu(D)&=\text{number of parts of $D$ for which }D_{ij}>j-i,\\[1mm]
\label{muD}\mu(D)&=\text{number of parts of $D$ for which }D_{ij}\le j-i,\\[1mm]
\label{rho1D}\rho_1(D)&=\text{number of $n$'s in $D$},\\[1mm]
\label{rho2D}\rho_2(D)&=\text{(number of $(n-1)$'s in $D$)} + \text{(number of rows of $D$ of length }n-1).
\end{align}

It follows from~\eqref{rowcond} and~\eqref{ASMDPPn} that, for the DPP statistic~$\rho_1(D)$,
$n$'s can appear in only the first row of~$D$, and
that, for the DPP statistic~$\rho_2(D)$,
$(n-1)$'s can appear in only the first two rows of $D$, and the number of rows of~$D$ of length $n-1$ is either~0 or~1,
since only the first row can have this length.

For an ASM~$A$, $\nu(A)$ (which can easily be shown to be nonnegative) is regarded as the number of generalized inversions in~$A$, since
if $A$ is a permutation matrix, i.e., if $\mu(A)=0$, then $\nu(A)$ is the number of inversions
in the permutation $\pi$ given by $\delta_{\pi_i,j}=A_{ij}$.
This statistic can also be written as
$\nu(A)=\sum_{1\le i\le i'\le n;\;1\le j'<j\le n}A_{ij}\,A_{i'j'}$,
where this can be derived using the fact that the sum of entries in each row and column of $A$ is a constant.

The definitions of the statistics $\nu$, $\mu$ and $\rho_1$ for ASMs and DPPs are based on definitions of
statistics introduced by Mills, Robbins and Rumsey~\cite[pp.~344--345]{MilRobRum83}.
There are some minor differences between
the definitions used here and those of Mills, Robbins and Rumsey, with full details of these differences given in~\cite[Sec.~1.2]{BehDifZin12a}.
For example, for an ASM or DPP~$X$, Mills Robbins and Rumsey use $\nu(X)+\mu(X)$ instead of $\nu(X)$ as one of the basic statistics.
It can be seen that, for an ASM~$A$, $\nu(A)+\mu(A)$ can be regarded as an alternative generalization of the number of inversions in
a permutation, since~$\mu(A)=0$ if $A$ is a permutation matrix.

For a DPP~$D$, the parts $D_{ij}$ which contribute to $\nu(D)$ and $\mu(D)$ in~\eqref{nuD}--\eqref{muD} are
referred to as the nonspecial and special parts of~$D$ respectively,
where this terminology matches that introduced by Mills, Robbins and Rumsey~\cite[p.~344]{MilRobRum83}.
It can be seen that a DPP~$D$ contains at least as many nonspecial parts as rows, since the
first part of each row is nonspecial.

The statistic~$\rho_2$ for ASMs was
first considered, in the context of its joint distribution with~$\rho_1$, and a conjectured connection
with totally symmetric self-complementary plane partitions, by Mills, Robbins and Rumsey~\cite[p.~284]{MilRobRum86}.
The statistic~$\rho_2$ for DPPs has not previously been considered.

It can be seen that, for an ASM or DPP~$X$, $\nu(X)$ and $\mu(X)$ depend on the bulk structure of~$X$,
while $\rho_1(X)$ and $\rho_2(X)$ depend only on the structure at or near the boundary of~$X$.

The statistics for the examples in \eqref{ASMDPPEx} are, taking $n=6$ for the DPP $D$,
\begin{alignat}{4}\notag\nu(A)&=5,\;\;\;&\mu(A)&=3,\;\;\;&\rho_1(A)&=3,\;\;\;&\rho_2(A)&=2,\\
\label{ASMDPPstatex}\nu(D)&=7,\;\;&\mu(D)&=2,\;\;&\rho_1(D)&=3,\;\;&\rho_2(D)&=2.\end{alignat}

Doubly-refined (or four-statistic) generating functions, which give weighted enumerations of
the elements of $\ASM(n)$ or~$\DPP(n)$ using arbitrary weights~$x$,~$y$,~$z_1$ and~$z_2$
associated with the statistics~\eqref{nuA}--\eqref{rho2D}, are now defined as
\begin{align}\label{ZASM}Z^\ASM_n(x,y,z_1,z_2)&=\sum_{A\in\ASM(n)}x^{\nu(A)}\,y^{\mu(A)}\,z_1^{\rho_1(A)}\,z_2^{\rho_2(A)},\\[2mm]
\label{ZDPP}Z^\DPP_n(x,y,z_1,z_2)&=\sum_{D\in\DPP(n)}x^{\nu(D)}\,y^{\mu(D)}\,z_1^{\rho_1(D)}\,z_2^{\rho_2(D)}.\end{align}
As indicated in Section~\ref{intro}, the term doubly-refined refers to the fact that these generating functions each involve two
boundary statistics.

As examples of~\eqref{ZASM}--\eqref{ZDPP}, \eqref{ASMDPP3} gives
\begin{multline}\label{ZASMDPP3}Z^\ASM_3(x,y,z_1,z_2)=Z^\DPP_3(x,y,z_1,z_2)=\\1+x^3z_1^2z_2^2+xz_2+x^2z_1^2z_2+xz_1+x^2z_1z_2^2+xyz_1z_2,\end{multline}
where the terms are written in an order which corresponds to that used in each set in~\eqref{ASMDPP3}.

Singly-refined (or three-statistic) generating functions for ASMs and DPPs are now defined as
\begin{equation}\label{Zsing}Z^\ASM_n(x,y,z)=Z^\ASM_n(x,y,z,1),\qquad Z^\DPP_n(x,y,z)=Z^\DPP_n(x,y,z,1).\end{equation}

The singly- and doubly-refined generating functions are also related, for $n\ge2$, by
\begin{equation}\label{Z0}Z^\ASM_{n-1}(x,y,z)=Z^\ASM_n(x,y,z,0),\qquad Z^\DPP_{n-1}(x,y,z)=Z^\DPP_n(x,y,z,0).\end{equation}
These identities can be proved by constructing simple bijections between
$\{A\in\ASM(n-1)\mid\nu(A)=p,\;\mu(A)=m,\;\rho_1(A)=k\}$ and
$\{A\in\ASM(n)\mid\nu(A)=p,\;\mu(A)=m,\;\rho_1(A)=k,\;\rho_2(A)=0\}$, and between
$\{D\in\DPP(n-1)\mid\nu(D)=p,\;\mu(D)=m,\;\rho_1(D)=k\}$ and
$\{D\in\DPP(n)\mid\nu(D)=p,\;\mu(D)=m,\;\rho_1(D)=k,\;\rho_2(D)=0\}$.  Specifically,
an ASM in the first set is mapped to the second set by augmentation with a row and column on the bottom and right,
where the common entry of the additional row and column is~$1$ and all other entries are~$0$,
an ASM in the second set is mapped to the first set by deletion of its last row and last column,
a DPP in the first set is mapped to the second set by replacing each of its $(n-1)$'s by an $n$,
and a DPP in the second set is mapped to the first set by replacing each of its $n$'s by an $n-1$.

It can also be shown that $Z^\ASM_{n-1}(x,y,z)=Z^\ASM_n(x,y,0,z)$ and $Z^\DPP_{n-1}(x,y,z)=Z^\DPP_n(x,y,$
$0,z)$, so that $Z^\ASM_{n-2}(x,y,1)=Z^\ASM_n(x,y,0,0)$ and
$Z^\DPP_{n-2}(x,y,1)=Z^\DPP_n(x,y,0,0)$, although these relations are not needed for the main
proofs of this paper.

\section{Main results}\label{mainres}
In this section, the main results of this paper are stated, and their connections to the main results of the
paper~\cite{BehDifZin12a} are outlined.  The overall method of proving the results is also described, but
the presentation of the details of these proofs is deferred to later sections.

It was conjectured by Mills, Robbins and Rumsey~\cite[Conj.~3]{MilRobRum83}, and
proved in~\cite{BehDifZin12a}, that the ASM and DPP singly-refined generating functions are equal.
\begin{theorem*}[See \protect{\cite[Thm.~1]{BehDifZin12a}}]
For any $n$, $x$, $y$ and $z$,
\begin{equation}\label{MRR}Z^\ASM_n(x,y,z)=Z^\DPP_n(x,y,z).\end{equation}
\end{theorem*}
A subsidiary result of~\cite[Eqs.~(97)--(98)]{BehDifZin12a}
is that the singly-refined generating functions can be expressed as the determinant of a certain $n\times n$ matrix. Specifically,
\begin{equation}\label{singdet}Z^\ASM_n(x,y,z)=Z^\DPP_n(x,y,z)=\det_{0\le i,j\le n-1}\bigl(K_n(x,y,z)_{ij}\bigr),\end{equation}
where
\begin{equation}\label{Ksing}K_n(x,y,z)_{ij}=-\delta_{i,j+1}+
\begin{cases}
\sum_{k=0}^{\min(i,j+1)}\binom{i-1}{i-k}\binom{j+1}{k}x^ky^{i-k},&j\le n-2,\\[1.5mm]
\sum_{k=0}^i\sum_{l=0}^k\binom{i-1}{i-k}\binom{n-l-1}{k-l}x^ky^{i-k}z^l,&j=n-1.
\end{cases}
\end{equation}
(For other, transformed versions of this determinant formula, based on formulae for the enumeration of DPPs obtained by
Mills, Robbins and Rumsey~\cite[p.~346]{MilRobRum83}
and Lalonde~\cite[Thm.~3.1]{Lal02}, see, for example,~\cite[Eqs.~(65)--(66) \&~(87)--(88)]{BehDifZin12a}.)

The primary result of this paper is that the ASM and DPP doubly-refined generating functions are also equal.
\begin{theorem}\label{DRT}For any $n$, $x$, $y$, $z_1$ and $z_2$,
\begin{equation}\label{DRT1}Z^\ASM_n(x,y,z_1,z_2)=Z^\DPP_n(x,y,z_1,z_2).\end{equation}
Equivalently, for any~$n$,~$p$,~$m$,~$k_1$ and $k_2$, the sizes of
$\{A\in\ASM(n)\mid\nu(A)=p,\;\mu(A)=m,\;\rho_1(A)=k_1,\;\rho_2(A)=k_2\}$
and
$\{D\in\DPP(n)\mid\nu(D)=p,\;\mu(D)=m,\;\rho_1(D)=k_1,\;\rho_2(D)=k_2\}$
are equal.\end{theorem}
It can be seen, using the definitions of ASMs and DPPs and the statistics~\eqref{nuA}--\eqref{rho2D}, that each set
in Theorem~\ref{DRT} is empty unless the integers $p$, $m$, $k_1$ and $k_2$ lie within the ranges
\begin{equation}p=0,\ldots,\tfrac{n(n-1)}{2},\qquad m=0,\ldots,\begin{cases}\tfrac{(n-1)^2}{4},&n\text{ odd},\\
\tfrac{n(n-2)}{4},&n\text{ even,}\end{cases}\qquad k_1,k_2=0,\ldots,n-1.\end{equation}

A subsidiary result of this paper is that
the doubly-refined generating functions can be expressed as the determinant of a certain $n\times n$ matrix.  Specifically,
\begin{equation}\label{doubdet}Z^\ASM_n(x,y,z_1,z_2)=Z^\DPP_n(x,y,z_1,z_2)=\det_{0\le i,j\le n-1}\bigl(K_n(x,y,z_1,z_2)_{ij}\bigr),\end{equation}
where
\begin{multline}\label{K}K_n(x,y,z_1,z_2)_{ij}=\\
-\delta_{i,j+1}+\begin{cases}
\sum_{k=0}^{\min(i,j+1)}\binom{i-1}{i-k}\binom{j+1}{k}x^ky^{i-k},&j\le n-3,\\[1.5mm]
\sum_{k=0}^i\sum_{l=0}^k\binom{i-1}{i-k}\binom{n-l-2}{k-l}x^ky^{i-k}z_2^{l+1},&j=n-2,\\[1.5mm]
\sum_{k=0}^i\sum_{l=0}^k\sum_{m=0}^l\binom{i-1}{i-k}\binom{n-l-2}{k-l}x^ky^{i-k}z_1^mz_2^{l-m},&j=n-1.
\end{cases}\end{multline}
(For an alternative, transformed version of this determinant formula see~\eqref{Ldet}.)
The matrices of~\eqref{Ksing} and~\eqref{K} are related by
$K_n(x,y,z)=K_n(x,y,z,1)$.

Theorem~\ref{DRT} and~\eqref{doubdet} are valid for the trivial case $n=1$ (for which $\ASM(1)=\{(1)\}$, $\DPP(1)=\{\emptyset\}$
and $Z^\ASM_1(x,y,z_1,z_2)=Z^\DPP_1(x,y,z_1,z_2)=1$).  However, this case would require some distracting qualifications to be made to the statement
of some subsequent results. It will therefore be assumed, throughout the remainder of this paper, that $n\ge2$.

Theorem~\ref{DRT} will be proved by deriving identities
which express the ASM and DPP doubly-refined generating functions as identical
combinations of their respective singly-refined generating functions.
The required equality~\eqref{DRT1} of the doubly-refined generating functions
then follows immediately from the already-known equality~\eqref{MRR} of the singly-refined generating functions.
The identities expressing the doubly-refined generating functions in terms of their singly-refined counterparts
are given by~\eqref{propeq1} in the following theorem.
\begin{theorem}\label{prop}Let $Z$ denote $Z^\ASM$ or $Z^\DPP$.
Then, for any $n$, $x$, $y$, $z_1$ and $z_2$,
\begin{multline}\label{propeq1}
(z_1-z_2)\,Z_n(x,y,z_1,z_2)\,Z_{n-1}(x,y,1)=(z_1-1)z_2\,Z_n(x,y,z_1)\,Z_{n-1}(x,y,z_2)\;-\\
z_1(z_2-1)\,Z_{n-1}(x,y,z_1)\,Z_n(x,y,z_2).
\end{multline}
Equivalently, for any $n$, $x$, $y$, $z_1$, $z_2$, $z_3$ and $z_4$,
\begin{multline}\label{propeq2}
(z_1-z_2)\,(z_3-z_4)\,Z_n(x,y,z_1,z_2)\,Z_n(x,y,z_3,z_4)\;-\\
(z_1-z_3)\,(z_2-z_4)\,Z_n(x,y,z_1,z_3)\,Z_n(x,y,z_2,z_4)\;+\\
(z_1-z_4)\,(z_2-z_3)\,Z_n(x,y,z_1,z_4)\,Z_n(x,y,z_2,z_3)=0.
\end{multline}
\end{theorem}
The ASM case of~\eqref{propeq1} with $x=y=1$ was first obtained by Stroganov~\cite[Sec.~5]{Str06},
and the general ASM case of~\eqref{propeq1} was first obtained by Colomo and
Pronko~\cite[Eq.~(5.32)]{ColPro05b},~\cite[Eq.~(3.32)]{ColPro06}.  In the latter
case,~\eqref{propeq1} is expressed using the terminology of one- and two-point boundary
correlation functions for the six-vertex model with domain-wall boundary
conditions.  Results for such boundary correlation functions have also been obtained by Foda and Preston~\cite{FodPre04}.

The equivalence of \eqref{propeq1} and~\eqref{propeq2} can be verified as follows.
Identity~\eqref{propeq2} follows from~\eqref{propeq1} by
multiplying the LHS of~\eqref{propeq2} by $Z_{n-1}(x,y,1)^2$, using
\eqref{propeq1} to express each of the six cases of $(z_i-z_j)Z_n(x,y,z_i,z_j)Z_{n-1}(x,y,1)$
in terms of singly-refined generating functions,
and then checking that the resulting expression vanishes.
Conversely,~\eqref{propeq1} follows from~\eqref{propeq2}
by setting $z_3=1$ and $z_4=0$, and then using~\eqref{Zsing}--\eqref{Z0}.

It also follows from~\eqref{propeq1} that
\begin{equation}\label{ASMDPPint}Z^\ASM_n(x,y,z_1,z_2)=Z^\ASM_n(x,y,z_2,z_1),\qquad
Z^\DPP_n(x,y,z_1,z_2)=Z^\DPP_n(x,y,z_2,z_1).\end{equation}
This symmetry in $z_1$ and $z_2$ can alternatively be derived using simple operations on ASMs and DPPs, as will be done in Section~\ref{symm}.

Proofs of the ASM and DPP cases of Theorem~\ref{prop} will be given in Sections~\ref{ASMpropsect} and~\ref{DPPpropsect} respectively.
Each of these proofs will use a form of the Desnanot--Jacobi determinant identity, which will be stated in Section~\ref{DesJacSec}.
This identity will be applied to determinant formulae
which will be obtained using standard methods, as also used in~\cite{BehDifZin12a}
for parts of the proof of the equality~\eqref{MRR} of the singly-refined generating functions.
Specifically, for the ASM case of Theorem~\ref{prop},
a bijection between ASMs and configurations of the six-vertex model with domain-wall boundary conditions,
together with the Izergin--Korepin determinant formula for the partition function of this model, will be used, and,
for the DPP case of Theorem~\ref{prop}, a bijection between DPPs and certain sets of nonintersecting lattice paths,
together with the Lindstr\"{o}m--Gessel--Viennot theorem for the weighted enumeration of nonintersecting paths,
will be used.  The resulting determinant formula for the DPP case, which will be obtained in
Section~\ref{ZDPPsect}, will provide the second equality of~\eqref{doubdet}.

The proof given here of the ASM case of Theorem~\ref{prop} will differ from the proof given by Colomo and
Pronko~\cite{ColPro05b,ColPro06}.

\section{The Desnanot--Jacobi identity}\label{DesJacSec}
In this section, the determinant identity, often referred to as the Desnanot--Jacobi identity, is
stated, together with some background information.
This identity will then be used in Sections~\ref{ASMpropsect} and~\ref{DPPpropsect}.

For a matrix~$M$, and subsets $I$ and $J$ of the index sets for the rows and columns of~$M$,
let~$M^I_J$ denote the submatrix of $M$ in which the
rows of $I$ and columns of $J$ have been deleted.  Omission of $I$ or $J$ in this notation means that the
corresponding set is empty.
\begin{theorem*}[Desnanot, Jacobi]
For any $n\times n$ matrix $\bigl(M_{ij}\bigr)_{1\le i,j\le n}$, $1\le i_1<i_2\le n$ and $1\le j_1<j_2\le n$,
\begin{equation}\label{DesJac}\det M\,\det M^{\{i_1,i_2\}}_{\{j_1,j_2\}}=\det M^{\{i_1\}}_{\{j_1\}}\,\det M^{\{i_2\}}_{\{j_2\}}-
\det M^{\{i_1\}}_{\{j_2\}}\,\det M^{\{i_2\}}_{\{j_1\}}.
\end{equation}
Equivalently, for any $(n+2)\times n$ matrix $\bigl(N_{ij}\bigr)_{1\le i\le n+2;\,1\le j\le n}$ and
$1\le k_1<k_2<k_3<k_4\le n+2$,
\begin{equation}\label{Bazin}
\det N^{\{k_1,k_2\}}\,\det N^{\{k_3,k_4\}}-\det N^{\{k_1,k_3\}}\,\det N^{\{k_2,k_4\}}+\det N^{\{k_1,k_4\}}\,\det N^{\{k_2,k_3\}}=0.
\end{equation}
\end{theorem*}
In fact, only~\eqref{Bazin} will be directly used in this paper, but~\eqref{DesJac} is included here since
it represents (usually with $i_1=j_1=1$ and $i_2=j_2=2$ or $i_2=j_2=n$) the more commonly-stated form of the theorem.

The equivalence of~\eqref{DesJac} and~\eqref{Bazin} can be verified as follows.
Identity~\eqref{Bazin} can be obtained from~\eqref{DesJac} by augmenting the
$(n+2)\times n$ matrix $N$ with two columns on the left,
where these additional columns contain~$1$'s in row~$k_1$ of column~$1$ and row~$k_2$ of column~$2$
and~$0$'s elsewhere,
and then applying~\eqref{DesJac} with $i_1=k_3$, $i_2=k_4$, $j_1=1$ and $j_2=2$ to the $(n+2)\times(n+2)$ matrix.
Conversely,~\eqref{DesJac} can be
obtained from~\eqref{Bazin} by augmenting the $n\times n$ matrix $M$
with two rows on the top, where these additional rows contain~$1$'s in column~$j_1$ of row~1
and column~$j_2$ of row~$2$ and~$0$'s elsewhere,
and then applying~\eqref{Bazin} with $k_1=1$, $k_2=2$, $k_3=i_1+2$ and $k_4=i_2+2$
to the $(n+2)\times n$ matrix.

It can be shown similarly that~\eqref{DesJac} and~\eqref{Bazin} are also
equivalent to the identity that, for any $(n+1)\times n$ matrix $\bigl(P_{ij}\bigr)_{1\le i\le n+1;\,1\le j\le n}$,
$1\le k_1<k_2<k_3\le n+1$ and $1\le l\le n$,
\begin{equation}\label{DesJacBazin}
\det P^{\{k_1\}}\,\det P^{\{k_2,k_3\}}_{\{l\}}-\det P^{\{k_2\}}\,\det P^{\{k_1,k_3\}}_{\{l\}}+\det P^{\{k_3\}}\,\det P^{\{k_1,k_2\}}_{\{l\}}=0.
\end{equation}

An algebraic proof of~\eqref{DesJac} is given by
Bressoud~\cite[Sec.~3.5]{Bre99}, and a combinatorial proof of~\eqref{DesJac} is given by Zeilberger~\cite{Zei97}.

Cases of each of~\eqref{DesJac}--\eqref{DesJacBazin} for small values of~$n$ were published by Desnanot in~1819
(see Muir~\cite[Eqs.~(A)--(G), (A$'$)--(G$'$), pp.~139--142]{Mui06}).  The further attribution to Jacobi is based on the fact that~\eqref{DesJac}
for general~$n$ corresponds to the case $m=2$ of the identity,
published by Jacobi in~1834 (see Muir~\cite[Eq.~(XX.~4), p.~208]{Mui06}), that for any $n\times n$ matrix~$M$ and any $m\le n$, each
size $m$ minor of the matrix of size $n-1$ minors of~$M$ equals the complementary minor of~$M$ multiplied by $(\det M)^{m-1}$.
Proofs of the Jacobi identity, using
a variety of methods, are given, for example, by Brualdi and Schneider~\cite[Sec.~4]{BruSch83}, Knuth~\cite[Eq.~(3.16)]{Knu96},
Leclerc~\cite[Sec.~3.2]{Lec93}, Muir~\cite[Sec.~175]{Mui60} and Turnbull~\cite[pp.~77--79]{Tur60}.

In fact,~\eqref{DesJac}--\eqref{DesJacBazin} also correspond to special cases of various other determinant identities. For
example,~\eqref{DesJac} is a special case of an identity of Sylvester (see Brualdi and Berliner~\cite{BruBer08},
Brualdi and Schneider~\cite[Eq.~(8)]{BruSch83},
Knuth~\cite[Eq.~(3.17)]{Knu96} or Leclerc~\cite[Eq.~(8)]{Lec93}),
and~\eqref{Bazin} is a special case of an identity of Bazin (see Leclerc~\cite[Eq.~(7)]{Lec93}),
which is itself a special case of an identity of Reiss and Picquet (see Leclerc~\cite[Sec.~5.4]{Lec93}).
Furthermore,~\eqref{Bazin} can be regarded as a simple case of a Pl\"ucker relation (see, for example,
Harris~\cite[pp.~65--66]{Har92}).

It is interesting to note that the Desnanot--Jacobi identity is also related to ASMs
through a modified version of Dodgson's condensation algorithm.
In the standard algorithm~\cite{Dod66},
the determinant of an $n\times n$ matrix~$M$ is computed by successively taking $k=2,\ldots,n$, and using~\eqref{DesJac}
(with $i_1=j_1=1$ and $i_2=j_2=k$) to compute all connected minors of~$M$ of size~$k$, from the connected minors
of~$M$ of sizes~$k-1$ and~$k-2$ (with size~0 minors taken to have value~1). If the coefficient of
the second product of determinants on the RHS of~\eqref{DesJac} is changed to a parameter~$\lambda$, and this modified formula
is instead used throughout the algorithm, then the so-called $\lambda$-determinant of~$M$ is produced, and,
as shown by Robbins and Rumsey~\cite[Eq.~(27)]{RobRum86},
this is naturally expressed as a sum over $\ASM(n)$, rather than as a sum over permutations of $\{1,\ldots,n\}$.
For further details see also, for example, Abeles~\cite{Abe08}, Bressoud~\cite[Sec.~3.5]{Bre99} or Bressoud and Propp~\cite{BrePro99}.

\section{Proof of the ASM case of Theorem 2}\label{ASMpropsect}
\subsection{Bijection between ASMs and configurations of the six-vertex model with DWBC}\label{ASMbijsect}
In this section, the set of configurations of the statistical mechanical six-vertex
model on an $n\times n$ grid with domain-wall boundary conditions (DWBC) is described, and the
details of a natural bijection between $\ASM(n)$ and this set are summarized.  This is standard material, for which more information and
references can be found in, for example,~\cite[Secs.~2.1 \&~3.1]{BehDifZin12a}.

\begin{figure}[h]
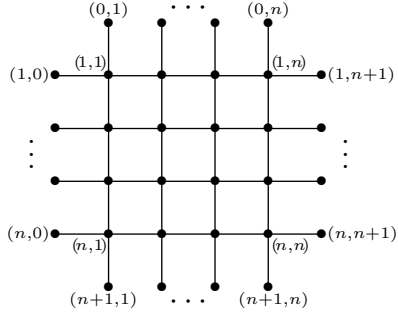
\centering\psset{unit=7mm}\pspicture(-0.9,-0.1)(6.4,5.5)
\multips(0,1)(0,1){4}{\psline[linewidth=0.5pt](0,0)(5,0)}\multips(1,0)(1,0){4}{\psline[linewidth=0.5pt](0,0)(0,5)}
\multirput(1,0)(1,0){4}{$\scriptstyle\bullet$}\multirput(0,1)(1,0){6}{$\scriptstyle\bullet$}
\multirput(0,2)(1,0){6}{$\scriptstyle\bullet$}\multirput(0,3)(1,0){6}{$\scriptstyle\bullet$}
\multirput(0,4)(1,0){6}{$\scriptstyle\bullet$}\multirput(1,5)(1,0){4}{$\scriptstyle\bullet$}
\rput[b](1,5.1){$\scriptscriptstyle(0,1)$}\rput[b](4,5.1){$\scriptscriptstyle(0,n)$}
\rput[t](0.9,-0.1){$\scriptscriptstyle(n+1,1)$}\rput[t](4.1,-0.1){$\scriptscriptstyle(n+1,n)$}
\rput[r](-0.1,4){$\scriptscriptstyle(1,0)$}\rput[r](-0.1,1){$\scriptscriptstyle(n,0)$}
\rput[l](5.1,1){$\scriptscriptstyle(n,n+1)$}\rput[l](5.1,4){$\scriptscriptstyle(1,n+1)$}
\rput[br](0.99,4.09){$\scriptscriptstyle(\!1,1\!)$}\rput[bl](4.05,4.09){$\scriptscriptstyle(\!1,n\!)$}
\rput[tl](4.05,0.91){$\scriptscriptstyle(\!n,n\!)$}\rput[tr](0.99,0.91){$\scriptscriptstyle(\!n,1\!)$}
\rput(2.54,5.28){$\cdots$}\rput(-0.45,2.69){$\vdots$}\rput(2.54,-0.28){$\cdots$}\rput(5.45,2.69){$\vdots$}\endpspicture
\caption{The undirected grid $\mathcal{G}_n$.}\label{grid1}\end{figure}

Let $\mathcal{G}_n$ be the $n\times n$ undirected grid with vertex set
$\{(i,j)\mid i,j=0,\ldots,n+1\}\setminus\{(0,0),(0,n+1),(n+1,0),(n+1,n+1)\}$,
where $(i,j)$ is taken to be in the $i$th row from the top and $j$th column from the left, and for which there are
horizontal edges between $(i,j)$ and $(i,j\pm1)$, and vertical edges between $(i,j)$ and $(i\pm1,j)$,
for each $i,j=1,\ldots,n$. This grid is shown in Figure~\ref{grid1}.
The descriptions `internal' and `external' are applied in the obvious way to the
vertices and edges of $\mathcal{G}_n$.

A configuration of the six-vertex model on $\mathcal{G}_n$ with DWBC is an assignment of arrows to the edges of
$\mathcal{G}_n$ in which the external edge arrows on the upper, right, lower and left boundaries of $\mathcal{G}_n$
are all directed upward, leftward, downward and rightward respectively, while the arrows on the four edges incident to any
internal vertex satisfy the condition that two point towards and two point away from the vertex.

Now define
\begin{equation}\SVDWBC(n)=\{\text{configurations of the six-vertex model on $\mathcal{G}_n$ with DWBC}\}.\end{equation}
For example,
\psset{unit=4.4mm}
\begin{multline}\label{6VDBWC3}\SVDWBC(3)=\\
\left\{\raisebox{-7.4mm}{
\pspicture(0.4,0.1)(4.5,3.9)\multips(0.1,1)(0,1){3}{\psline[linewidth=0.5pt](0,0)(3.8,0)}\multips(1,0.1)(1,0){3}{\psline[linewidth=0.5pt](0,0)(0,3.8)}
\psdots[dotstyle=triangle*,dotscale=1,dotangle=0](1,3.5)(2,2.5)(2,3.5)(3,1.5)(3,2.5)(3,3.5)
\psdots[dotstyle=triangle*,dotscale=1,dotangle=180](1,0.5)(1,1.5)(1,2.5)(2,0.5)(2,1.5)(3,0.5)
\psdots[dotstyle=triangle*,dotscale=1,dotangle=90](1.5,3)(2.5,2)(2.5,3)(3.5,1)(3.5,2)(3.5,3)
\psdots[dotstyle=triangle*,dotscale=1,dotangle=-90](0.5,1)(0.5,2)(0.5,3)(1.5,1)(1.5,2)(2.5,1)
\rput(4.1,1.2){,}\endpspicture
\pspicture(0.1,0.1)(4.5,3.9)\multips(0.1,1)(0,1){3}{\psline[linewidth=0.5pt](0,0)(3.8,0)}\multips(1,0.1)(1,0){3}{\psline[linewidth=0.5pt](0,0)(0,3.8)}
\psdots[dotstyle=triangle*,dotscale=1,dotangle=0](1,1.5)(1,2.5)(1,3.5)(2,2.5)(2,3.5)(3,3.5)
\psdots[dotstyle=triangle*,dotscale=1,dotangle=180](1,0.5)(2,0.5)(2,1.5)(3,0.5)(3,1.5)(3,2.5)
\psdots[dotstyle=triangle*,dotscale=1,dotangle=90](1.5,1)(2.5,1)(2.5,2)(3.5,1)(3.5,2)(3.5,3)
\psdots[dotstyle=triangle*,dotscale=1,dotangle=-90](0.5,1)(0.5,2)(0.5,3)(1.5,2)(1.5,3)(2.5,3)\rput(4.1,1.2){,}\endpspicture
\pspicture(0.1,0.1)(4.5,3.9)\multips(0.1,1)(0,1){3}{\psline[linewidth=0.5pt](0,0)(3.8,0)}\multips(1,0.1)(1,0){3}{\psline[linewidth=0.5pt](0,0)(0,3.8)}
\psdots[dotstyle=triangle*,dotscale=1,dotangle=0](1,3.5)(2,1.5)(2,2.5)(2,3.5)(3,2.5)(3,3.5)
\psdots[dotstyle=triangle*,dotscale=1,dotangle=180](1,0.5)(1,1.5)(1,2.5)(2,0.5)(3,0.5)(3,1.5)
\psdots[dotstyle=triangle*,dotscale=1,dotangle=90](1.5,3)(2.5,1)(2.5,3)(3.5,1)(3.5,2)(3.5,3)
\psdots[dotstyle=triangle*,dotscale=1,dotangle=-90](0.5,1)(0.5,2)(0.5,3)(1.5,1)(1.5,2)(2.5,2)\rput(4.1,1.2){,}\endpspicture
\pspicture(0.1,0.1)(4.5,3.9)\multips(0.1,1)(0,1){3}{\psline[linewidth=0.5pt](0,0)(3.8,0)}\multips(1,0.1)(1,0){3}{\psline[linewidth=0.5pt](0,0)(0,3.8)}
\psdots[dotstyle=triangle*,dotscale=1,dotangle=0](1,2.5)(1,3.5)(2,1.5)(2,2.5)(2,3.5)(3,3.5)
\psdots[dotstyle=triangle*,dotscale=1,dotangle=180](1,0.5)(1,1.5)(2,0.5)(3,0.5)(3,1.5)(3,2.5)
\psdots[dotstyle=triangle*,dotscale=1,dotangle=90](1.5,2)(2.5,1)(2.5,2)(3.5,1)(3.5,2)(3.5,3)
\psdots[dotstyle=triangle*,dotscale=1,dotangle=-90](0.5,1)(0.5,2)(0.5,3)(1.5,1)(1.5,3)(2.5,3)\rput(4.1,1.2){,}\endpspicture
\pspicture(0.1,0.1)(4.5,3.9)\multips(0.1,1)(0,1){3}{\psline[linewidth=0.5pt](0,0)(3.8,0)}\multips(1,0.1)(1,0){3}{\psline[linewidth=0.5pt](0,0)(0,3.8)}
\psdots[dotstyle=triangle*,dotscale=1,dotangle=0](1,2.5)(1,3.5)(2,3.5)(3,1.5)(3,2.5)(3,3.5)
\psdots[dotstyle=triangle*,dotscale=1,dotangle=180](1,0.5)(1,1.5)(2,0.5)(2,1.5)(2,2.5)(3,0.5)
\psdots[dotstyle=triangle*,dotscale=1,dotangle=90](1.5,2)(2.5,2)(2.5,3)(3.5,1)(3.5,2)(3.5,3)
\psdots[dotstyle=triangle*,dotscale=1,dotangle=-90](0.5,1)(0.5,2)(0.5,3)(1.5,1)(1.5,3)(2.5,1)\rput(4.1,1.2){,}\endpspicture
\pspicture(0.1,0.1)(4.5,3.9)\multips(0.1,1)(0,1){3}{\psline[linewidth=0.5pt](0,0)(3.8,0)}\multips(1,0.1)(1,0){3}{\psline[linewidth=0.5pt](0,0)(0,3.8)}
\psdots[dotstyle=triangle*,dotscale=1,dotangle=0](1,1.5)(1,2.5)(1,3.5)(2,3.5)(3,2.5)(3,3.5)
\psdots[dotstyle=triangle*,dotscale=1,dotangle=180](1,0.5)(2,0.5)(2,1.5)(2,2.5)(3,0.5)(3,1.5)
\psdots[dotstyle=triangle*,dotscale=1,dotangle=90](1.5,1)(2.5,1)(2.5,3)(3.5,1)(3.5,2)(3.5,3)
\psdots[dotstyle=triangle*,dotscale=1,dotangle=-90](0.5,1)(0.5,2)(0.5,3)(1.5,2)(1.5,3)(2.5,2)\rput(4.1,1.2){,}\endpspicture
\pspicture(0.1,0.1)(4,3.9)\multips(0.1,1)(0,1){3}{\psline[linewidth=0.5pt](0,0)(3.8,0)}\multips(1,0.1)(1,0){3}{\psline[linewidth=0.5pt](0,0)(0,3.8)}
\psdots[dotstyle=triangle*,dotscale=1,dotangle=0](1,2.5)(1,3.5)(2,1.5)(2,3.5)(3,2.5)(3,3.5)
\psdots[dotstyle=triangle*,dotscale=1,dotangle=180](1,0.5)(1,1.5)(2,0.5)(2,2.5)(3,0.5)(3,1.5)
\psdots[dotstyle=triangle*,dotscale=1,dotangle=90](1.5,2)(2.5,1)(2.5,3)(3.5,1)(3.5,2)(3.5,3)
\psdots[dotstyle=triangle*,dotscale=1,dotangle=-90](0.5,1)(0.5,2)(0.5,3)(1.5,1)(1.5,3)(2.5,2)\endpspicture}
\right\}.\end{multline}

\begin{figure}[h]
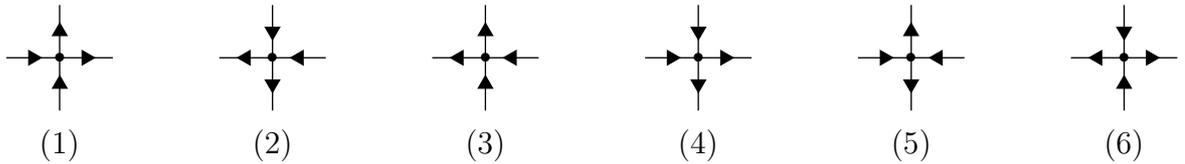
\centering\psset{unit=14mm}\pspicture(0.5,-0.3)(11.5,1)
\multips(1,0.5)(2,0){6}{\psline[linewidth=0.5pt](0,-0.5)(0,0.5)\psline[linewidth=0.5pt](-0.5,0)(0.5,0)}
\multirput(1,0.5)(2,0){6}{$\scriptstyle\bullet$}
\rput(1,-0.33){(1)}\psdots[dotstyle=triangle*,dotscale=1.5](1,0.75)\psdots[dotstyle=triangle*,dotscale=1.5](1,0.25)
\psdots[dotstyle=triangle*,dotscale=1.5,dotangle=-90](0.75,0.5)\psdots[dotstyle=triangle*,dotscale=1.5,dotangle=-90](1.25,0.5)
\rput(3,-0.33){(2)}\psdots[dotstyle=triangle*,dotscale=1.5,dotangle=180](3,0.75)\psdots[dotstyle=triangle*,dotscale=1.5,dotangle=180](3,0.25)
\psdots[dotstyle=triangle*,dotscale=1.5,dotangle=90](2.75,0.5)\psdots[dotstyle=triangle*,dotscale=1.5,dotangle=90](3.25,0.5)
\rput(5,-0.33){(3)}\psdots[dotstyle=triangle*,dotscale=1.5](5,0.75)\psdots[dotstyle=triangle*,dotscale=1.5](5,0.25)
\psdots[dotstyle=triangle*,dotscale=1.5,dotangle=90](4.75,0.5)\psdots[dotstyle=triangle*,dotscale=1.5,dotangle=90](5.25,0.5)
\rput(7,-0.33){(4)}\psdots[dotstyle=triangle*,dotscale=1.5,dotangle=180](7,0.75)\psdots[dotstyle=triangle*,dotscale=1.5,dotangle=180](7,0.25)
\psdots[dotstyle=triangle*,dotscale=1.5,dotangle=-90](6.75,0.5)\psdots[dotstyle=triangle*,dotscale=1.5,dotangle=-90](7.25,0.5)
\rput(9,-0.33){(5)}\psdots[dotstyle=triangle*,dotscale=1.5](9,0.75)\psdots[dotstyle=triangle*,dotscale=1.5,dotangle=180](9,0.25)
\psdots[dotstyle=triangle*,dotscale=1.5,dotangle=-90](8.75,0.5)\psdots[dotstyle=triangle*,dotscale=1.5,dotangle=90](9.25,0.5)
\rput(11,-0.33){(6)}\psdots[dotstyle=triangle*,dotscale=1.5,dotangle=180](11,0.75)\psdots[dotstyle=triangle*,dotscale=1.5](11,0.25)
\psdots[dotstyle=triangle*,dotscale=1.5,dotangle=90](10.75,0.5)\psdots[dotstyle=triangle*,dotscale=1.5,dotangle=-90](11.25,0.5)
\endpspicture\caption{The possible arrow configurations on edges incident to an internal vertex.}\label{vertices}\end{figure}

In an element of $\SVDWBC(n)$,
there are six possible configurations of arrows on the four edges incident to an internal vertex of $\mathcal{G}_n$.  These types of configurations,
known as vertex configurations, are shown in Figure~\ref{vertices}, where the numbers will be used to label the types.

For any $C\in\SVDWBC(n)$, denote the total number of type-$k$ vertex configurations in~$C$ as~$\mathcal{N}_{(k)}(C)$,
and denote the number of type-$k$ vertex configurations in~$C$ in row~$i$ of $\mathcal{G}_n$ as~$\mathcal{N}^i_{(k)}(C)$.

It can be shown that, for any $C\in\SVDWBC(n)$,
\begin{gather}\notag\mathcal{N}_{(1)}(C)=\mathcal{N}_{(2)}(C),\quad\mathcal{N}_{(3)}(C)=\mathcal{N}_{(4)}(C),\quad
\mathcal{N}_{(5)}(C)=\mathcal{N}_{(6)}(C)+n,\\
\notag\mathcal{N}^1_{(2)}(C)=\mathcal{N}^1_{(4)}(C)=\mathcal{N}^1_{(6)}(C)=
\mathcal{N}^n_{(1)}(C)=\mathcal{N}^n_{(3)}(C)=\mathcal{N}^n_{(6)}(C)=0,\\
\label{Nrel}\mathcal{N}^1_{(5)}(C)=\mathcal{N}^n_{(5)}(C)=1.\end{gather}

Statistics for each $C\in\SVDWBC(n)$ are now defined as
\begin{equation}\label{6VDWBCstat}\nu(C)=\mathcal{N}_{(1)}(C),\quad\mu(C)=\mathcal{N}_{(6)}(C),\quad
\rho_1(C)=\mathcal{N}^1_{(1)}(C),\quad\rho_2(C)=\mathcal{N}^n_{(2)}(C).\end{equation}

It can be shown straightforwardly that there is a natural bijection between $\ASM(n)$ and $\SVDWBC(n)$,
and that for each $A\in\ASM(n)$ and $C\in\SVDWBC(n)$ which correspond under this bijection,
$\nu(A)=\nu(C)$, $\mu(A)=\mu(C)$, $\rho_1(A)=\rho_1(C)$ and $\rho_2(A)=\rho_2(C)$
(where the statistics in these equations are given
by~\eqref{nuA}--\eqref{rho2A} on the LHS and~\eqref{6VDWBCstat} on the RHS).

The details of this bijection are as follows.
To map $A\in\ASM(n)$ to $C\in\SVDWBC(n)$, first associate the partial row sum $\sum_{j'=1}^{j}A_{ij'}$
with the horizontal edge between $(i,j)$ and $(i,j+1)$,
for each $i=1,\ldots,n$, $j=0,\ldots,n$, and associate
the partial column sum $\sum_{i'=1}^{i}A_{i'j}$ with the vertical edge between $(i,j)$ and $(i+1,j)$,
for each $i=0,\ldots,n$, $j=1,\ldots,n$.
(The defining properties of ASMs imply that each of these partial sums is~0 or~1.)  Then
obtain $C$ by assigning a rightward or upward arrow to each edge associated with a~$0$,
and a leftward or downward arrow to each edge associated with a~$1$.
To map $C\in\SVDWBC(n)$ to $A\in\ASM(n)$, set $A_{ij}$ to be~$1$,~$-1$ or~$0$ according to whether the vertex configuration in~$C$
at internal vertex $(i,j)$ is of type~5,~6 or~1--4 respectively.

For example, the ASM~$A$ of~\eqref{ASMDPPEx} and the corresponding six-vertex model configuration are
\begin{equation}\label{6VDWBCex}\left(\begin{array}{@{}c@{\;\;}c@{\;\;}c@{\;\;}c@{\;\;}c@{\;\;}c@{}}0&0&0&1&0&0\\[0.2mm]0&1&0&-1&1&0\\[0.2mm]1&-1&1&0&0&0
\\[0.2mm]0&0&0&1&0&0\\[0.2mm]0&1&0&-1&0&1\\[0.2mm]0&0&0&1&0&0\end{array}\right)\quad\longleftrightarrow\quad
\psset{unit=5.1mm}\raisebox{-18mm}{\pspicture(-0.1,-0.25)(7,7)
\multips(0.1,1)(0,1){6}{\psline[linewidth=0.5pt](0,0)(6.8,0)}\multips(1,0.1)(1,0){6}{\psline[linewidth=0.5pt](0,0)(0,6.8)}
\multirput(1,0)(1,0){6}{$\scriptscriptstyle\bullet$}\multirput(0,1)(1,0){8}{$\scriptscriptstyle\bullet$}
\multirput(0,2)(1,0){8}{$\scriptscriptstyle\bullet$}\multirput(0,3)(1,0){8}{$\scriptscriptstyle\bullet$}
\multirput(0,4)(1,0){8}{$\scriptscriptstyle\bullet$}\multirput(0,5)(1,0){8}{$\scriptscriptstyle\bullet$}
\multirput(0,6)(1,0){8}{$\scriptscriptstyle\bullet$}\multirput(1,7)(1,0){6}{$\scriptscriptstyle\bullet$}\rput(7.6,3.5){.}
\psdots[dotstyle=triangle*,dotscale=1.1,dotangle=0]
(1,6.5)(2,6.5)(3,6.5)(4,6.5)(5,6.5)(6,6.5)(1,5.5)(2,5.5)(3,5.5)(5,5.5)(6,5.5)
(1,4.5)(3,4.5)(4,4.5)(6,4.5)(2,3.5)(4,3.5)(6,3.5)(2,2.5)(6,2.5)(4,1.5)
\psdots[dotstyle=triangle*,dotscale=1.1,dotangle=180]
(1,0.5)(2,0.5)(3,0.5)(4,0.5)(5,0.5)(6,0.5)(1,1.5)(2,1.5)(3,1.5)(5,1.5)(6,1.5)
(1,2.5)(3,2.5)(4,2.5)(5,2.5)(1,3.5)(3,3.5)(5,3.5)(2,4.5)(5,4.5)(4,5.5)
\psdots[dotstyle=triangle*,dotscale=1.1,dotangle=90]
(6.5,1)(6.5,2)(6.5,3)(6.5,4)(6.5,5)(6.5,6)(5.5,1)(5.5,3)(5.5,4)(5.5,5)(5.5,6)
(4.5,1)(4.5,3)(4.5,4)(4.5,6)(3.5,2)(3.5,4)(3.5,5)(2.5,2)(2.5,5)(1.5,4)
\psdots[dotstyle=triangle*,dotscale=1.1,dotangle=-90]
(0.5,1)(0.5,2)(0.5,3)(0.5,4)(0.5,5)(0.5,6)(1.5,1)(1.5,2)(1.5,3)(1.5,5)(1.5,6)
(2.5,1)(2.5,3)(2.5,4)(2.5,6)(3.5,1)(3.5,3)(3.5,6)(4.5,2)(4.5,5)(5.5,2)\endpspicture}
\end{equation}
It can be seen that the statistics~\eqref{6VDWBCstat} for the
six-vertex model configuration $C$ in~\eqref{6VDWBCex} are $\nu(C)=5$, $\mu(C)=3$, $\rho_1(C)=3$ and $\rho_2(C)=2$, and that these match the
respective statistics~\eqref{nuA}--\eqref{rho2A} for the ASM~$A$ in~\eqref{6VDWBCex}, as given in~\eqref{ASMDPPstatex}.

As further examples, in~\eqref{ASMDPP3} and~\eqref{6VDBWC3} the elements of $\ASM(3)$ and $\SVDWBC(3)$
are listed in an order for which respective elements correspond under the bijection of this section.
It can be seen that $Z^\ASM_3(x,y,z_1,z_2)$, as given in~\eqref{ZASMDPP3},
could now be obtained using this bijection and~\eqref{6VDBWC3}.

\subsection{The partition function of the six-vertex model with DWBC}
In this section, the partition function of the six-vertex model with DWBC is introduced.  A relation between this
partition function, at certain values of its parameters, and the doubly-refined
ASM generating function~\eqref{ZASM} is then derived using the bijection of Section~\ref{ASMbijsect}.

Let a weight $W_{(k)}(u,v)$ be associated with the vertex configuration of type $k$, where $u$ and~$v$ are so-called spectral parameters.

The partition function for the case of the six-vertex model
of relevance here depends on these weights, and on parameters~$u_i$ and~$v_j$
associated with row~$i$ and column~$j$ of $\mathcal{G}_n$, for each $1\le i,j\le n$.  Specifically, this partition function is defined as
\begin{equation}\label{ZSVDWBC}Z^\SV_n(u_1,\ldots,u_n;v_1,\ldots,v_n)=
\sum_{C\in\SVDWBC(n)}\;\prod_{i,j=1}^n\,W_{(C_{ij})}(u_i,v_j),\end{equation}
where $C_{ij}$ is the type of vertex configuration in $C$ at internal vertex $(i,j)$.

Let the weights now satisfy
\begin{gather}\notag W_{(1)}(u,v)=W_{(2)}(u,v)=a(u,v),\quad W_{(3)}(u,v)=W_{(4)}(u,v)=b(u,v),\\
\label{Wabc}W_{(5)}(u,v)=W_{(6)}(u,v)=c(u,v),\end{gather}
for particular functions $a(u,v)$, $b(u,v)$ and $c(u,v)$.

If the spectral parameters in~\eqref{ZSVDWBC} are given by
\begin{equation}\label{specpar}
u_2=\ldots=u_{n-1}=r,\quad u_1=s,\quad u_n=t,\quad v_1=\ldots=v_n=w,\end{equation}
for parameters $r$, $s$, $t$ and $w$, then
\begin{multline*}
Z^\SV_n(s,r,\ldots,r,t;w,\ldots,w)=\\
\sum_{C\in\SVDWBC(n)}\!\!\textstyle a(r,w)^{2\mathcal{N}_{(1)}(C)}\,
\bigl(\frac{a(s,w)}{a(r,w)}\bigr)^{\mathcal{N}_{(1)}^1(C)}\,
\bigl(\frac{a(t,w)}{a(r,w)}\bigr)^{\mathcal{N}_{(2)}^n(C)}\,b(r,w)^{2\mathcal{N}_{(3)}(C)}\;\times\\[-3mm]
\shoveright{\textstyle\bigl(\frac{b(s,w)}{b(r,w)}\bigr)^{\mathcal{N}_{(3)}^1(C)}\,
\bigl(\frac{b(t,w)}{b(r,w)}\bigr)^{\mathcal{N}_{(4)}^n(C)}\,
c(r,w)^{2\mathcal{N}_{(6)}(C)+n-2}\,c(s,w)\,c(t,w)}\\
=b(r,w)^{(n-1)(n-2)}\,(b(s,w)\,b(t,w))^{n-1}\,c(r,w)^{n-2}\,c(s,w)\,c(t,w)\;\times\\
\displaystyle\sum_{C\in\SVDWBC(n)}\textstyle\!
\bigl(\frac{a(r,w)}{b(r,w)}\bigr)^{2\mathcal{N}_{(1)}(C)}\,
\bigl(\frac{c(r,w)}{b(r,w)}\bigr)^{2\mathcal{N}_{(6)}(C)}\,
\bigl(\frac{a(s,w)\,b(r,w)}{a(r,w)\,b(s,w)}\bigr)^{\mathcal{N}_{(1)}^1(C)}\,
\bigl(\frac{a(t,w)\,b(r,w)}{a(r,w)\,b(t,w)}\bigr)^{\mathcal{N}_{(2)}^n(C)},\end{multline*}
where~\eqref{Nrel} and the facts that $\sum_{k=1}^6\mathcal{N}_{(k)}(C)=n^2$ and $\sum_{k=1}^6\mathcal{N}^i_{(k)}(C)=n$ were used.
The bijection of Section \ref{ASMbijsect}, \eqref{ZASM}  and \eqref{6VDWBCstat} now give
\begin{multline}\label{ZCZASM}
Z^\SV_n(s,r,\ldots,r,t;w,\ldots,w)=b(r,w)^{(n-1)(n-2)}\,(b(s,w)\,b(t,w))^{n-1}\,c(r,w)^{n-2}\,c(s,w)\,c(t,w)\\
\textstyle\times Z^\ASM_n\Bigl(\bigl(\frac{a(r,w)}{b(r,w)}\bigr)^2,
\bigl(\frac{c(r,w)}{b(r,w)}\bigr)^2,
\frac{a(s,w)\,b(r,w)}{a(r,w)\,b(s,w)},
\frac{a(t,w)\,b(r,w)}{a(r,w)\,b(t,w)}\Bigr).\end{multline}

\subsection{The Izergin--Korepin determinant formula}
In this section, the Izergin--Korepin formula
for the partition function~\eqref{ZSVDWBC}, with certain assignments of the weights~\eqref{Wabc}, is stated.

It was shown by Izergin~\cite{Ize87}, using certain results of Korepin~\cite{Kor82}, that if
the weights~\eqref{Wabc} satisfy the Yang--Baxter equation,
then the partition function~\eqref{ZSVDWBC} can be expressed in terms of the determinant of a certain $n\times n$ matrix.

Let the weights~\eqref{Wabc} be given by
\begin{equation}\label{wuv}\textstyle
a(u,v)=uq-\frac{v}{q},\qquad b(u,v)=\frac{u}{q}-vq,\qquad c(u,v)=\bigl(q^2-\frac{1}{q^2}\bigr)\,u^{1/2}\,v^{1/2},\end{equation}
where $q$ is a further, global parameter.  These weights can be shown to satisfy the Yang--Baxter equation, with the
resulting Izergin--Korepin determinant formula being given by the following result.
\begin{theorem*}[Izergin]
The partition function~\eqref{ZSVDWBC}, with weights~\eqref{Wabc} and~\eqref{wuv}, satisfies
\begin{equation}\label{Idet1}
Z^\SV_n(u_1,\ldots,u_n;v_1,\ldots,v_n)=
\frac{\prod_{i,j=1}^na(u_i,v_j)\,b(u_i,v_j)}{\prod_{1\le i<j\le n}(u_i-u_j)(v_j-v_i)}
\det_{1\le i,j\le n}\left(\frac{c(u_i,v_j)}{a(u_i,v_j)\,b(u_i,v_j)}\right).\end{equation}
\end{theorem*}
For a proof of~\eqref{Idet1} see, for example,
Bogoliubov, Pronko and Zvonarev~\cite[Sec.~4]{BogProZvo02},
or Izergin, Coker and Korepin~\cite[Sec.~5]{IzeCokKor92}.

It is apparent from~\eqref{Idet1} that $Z^\SV_n(u_1,\ldots,u_n;v_1,\ldots,v_n)$
is symmetric in $u_1,\ldots,u_n$ and in $v_1,\ldots,v_n$.
It can also be seen that, although the determinant and the denominator of the prefactor on the RHS of~\eqref{Idet1}
both vanish if $u_i=u_j$ or $v_i=v_j$ for some $i\ne j$, the RHS has a well-defined limit in these cases,
as a polynomial in $u_1^{1/2},\ldots,u_n^{1/2},v_1^{1/2},\ldots,v_n^{1/2}$.   Accordingly, it will be
valid to use~\eqref{Idet1} to derive properties of the partition function for the values~\eqref{specpar}.

\subsection{Application of the Desnanot--Jacobi identity}
In this section, the final steps in the proof of the ASM case of~\eqref{propeq2} are taken.  These involve using
the form~\eqref{Bazin} of the Desnanot--Jacobi identity, the Izergin--Korepin determinant formula~\eqref{Idet1},
and the relation~\eqref{ZCZASM} between the doubly-refined ASM generating function and
the partition function of the six-vertex model with DWBC.

Applying the form~\eqref{Bazin} of Desnanot--Jacobi identity
to the $(n+2)\times n$ matrix given by $\Bigl(\frac{c(u_i,v_j)}{a(u_i,v_j)\,b(u_i,v_j)}\Bigr)_{1\le i\le n+2;\,1\le j\le n}$,
and then applying the Izergin--Korepin determinant formula~\eqref{Idet1} to each of the six determinants which appear,
it follows that, for any $u_1,\ldots,u_{n+2}$, $v_1,\ldots,v_n$ and $1\le k_1<k_2<k_3<k_4\le n+2$,
\begin{multline}\label{ZBazin}
(u_{k_1}-u_{k_2})\,(u_{k_3}-u_{k_4})\,Z^\SV_n(u_1,\ldots,\widehat{u}_{k_1},\ldots,\widehat{u}_{k_2},\ldots,u_{n+2};v_1,\ldots,v_n)\;\times\\
Z^\SV_n(u_1,\ldots,\widehat{u}_{k_3},\ldots,\widehat{u}_{k_4},\ldots,u_{n+2};v_1,\ldots,v_n)\\
-(u_{k_1}-u_{k_3})\,(u_{k_2}-u_{k_4})\,Z^\SV_n(u_1,\ldots,\widehat{u}_{k_1},\ldots,\widehat{u}_{k_3},\ldots,u_{n+2};v_1,\ldots,v_n)\;\times\\
Z^\SV_n(u_1,\ldots,\widehat{u}_{k_2},\ldots,\widehat{u}_{k_4},\ldots,u_{n+2};v_1,\ldots,v_n)\\
+(u_{k_1}-u_{k_4})\,(u_{k_2}-u_{k_3})\,Z^\SV_n(u_1,\ldots,\widehat{u}_{k_1},\ldots,\widehat{u}_{k_4},\ldots,u_{n+2};v_1,\ldots,v_n)\;\times\\
Z^\SV_n(u_1,\ldots,\widehat{u}_{k_2},\ldots,\widehat{u}_{k_3},\ldots,u_{n+2};v_1,\ldots,v_n)=0,
\end{multline}
where $u_1,\ldots,\widehat{u}_i,\ldots,\widehat{u}_j,\ldots,u_{n+2}$ denotes the omission of $u_i$ and $u_j$ from
$u_1,\ldots,u_{n+2}$.
Note that~\eqref{ZBazin} would be satisfied by any function which has the form of~\eqref{Idet1}, for arbitrary functions $a(u,v)$,
$b(u,v)$ and $c(u,v)$.

Setting $u_5=\ldots=u_{n+2}=r$, $v_1=\ldots=v_n=w$ and $k_i=i$ in~\eqref{ZBazin} gives
\begin{multline}\label{ZSVDWBCdoubref}
(u_1-u_2)\,(u_3-u_4)\,Z^\SV_n(u_1,r,\ldots,r,u_2;w,\ldots,w)\,Z^\SV_n(u_3,r,\ldots,r,u_4;w,\ldots,w)\;-\\
(u_1-u_3)\,(u_2-u_4)\,Z^\SV_n(u_1,r,\ldots,r,u_3;w,\ldots,w)\,Z^\SV_n(u_2,r,\ldots,r,u_4;w,\ldots,w)\;+\\
(u_1-u_4)\,(u_2-u_3)\,Z^\SV_n(u_1,r,\ldots,r,u_4;w,\ldots,w)\,Z^\SV_n(u_2,r,\ldots,r,u_3;w,\ldots,w)=0,
\end{multline}
where the symmetry of $Z^\SV_n(u_1,\ldots,u_n;v_1,\ldots,v_n)$ in $u_1,\ldots,u_n$ has been used.

Now let parameters $x$, $y$ and $z_1,\ldots,z_4$ be given, using the functions of~\eqref{wuv}, as
\begin{equation}\label{xyz}\textstyle x=\bigl(\frac{a(r,w)}{b(r,w)}\bigr)^2,\qquad y=\bigl(\frac{c(r,w)}{b(r,w)}\bigr)^2,\qquad
z_i=\frac{a(u_i,w)\,b(r,w)}{a(r,w)\,b(u_i,w)},\quad i=1,\ldots,4.\end{equation}
It can be checked that this allows arbitrary $x$, $y$ and $z_1,\ldots,z_4$ to be expressed in terms
of parameters~$q$,~$r$,~$w$ and $u_1,\ldots,u_4$.
It can also be checked that, for $i,j=1,\ldots,4$,
\begin{equation}\label{uz}u_i-u_j=\textstyle\frac{a(r,w)\,b(u_i,w)\,b(u_j,w)}{b(r,w)\,(q^{-2}-q^2)\,w}\,(z_i-z_j).\end{equation}

The ASM case of~\eqref{propeq2}
now follows from~\eqref{ZSVDWBCdoubref}, using~\eqref{ZCZASM},~\eqref{uz}, and the fact that the arbitrary~$x$,~$y$
and $z_1,\ldots,z_4$ in~\eqref{propeq2} can be parameterized by~\eqref{xyz}.

It is interesting to note that the Desnanot--Jacobi identity was also used together with
an Izergin--Korepin determinant formula by
Korepin and Zinn-Justin~\cite[Sec.~3]{KorZin00} (see also Sogo~\cite[Sec.~4]{Sog93}) in the computation of
the thermodynamic limit of the six-vertex model with DWBC in its ferroelectric and disordered regimes.

\section{Proof of the DPP case of Theorem 2}\label{DPPpropsect}
\subsection{Bijection between DPPs and sets of nonintersecting lattice paths}\label{DPPbijsect}
In this section, the details of a bijection between $\DPP(n)$ and a set of certain sets of nonintersecting lattice
paths on an~$n\times n$ grid are summarized.  More information, background material and
references can be found in~\cite[Secs.~2.2 \&~3.2]{BehDifZin12a}.

\begin{figure}[h]
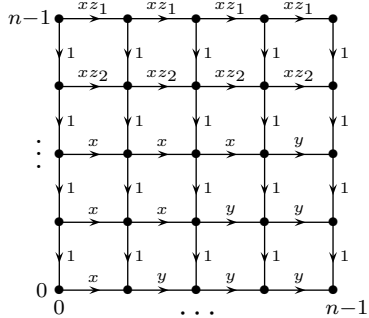
\centering
\psset{unit=4.5mm}
\pspicture(0,-0.4)(8.5,8.8)
\multips(0,0)(0,2){5}{\psline[arrows=->,arrowsize=3pt,linewidth=0.5pt](0,0)(1.2,0)
\multips(1,0)(2,0){3}{\psline[arrows=->,arrowsize=3pt,linewidth=0.5pt](0,0)(2.2,0)}\psline[linewidth=0.5pt](7,0)(8,0)}
\multips(0,8)(2,0){5}{\psline[arrows=->,arrowsize=3pt,linewidth=0.5pt](0,0)(0,-1.2)
\multips(0,-1)(0,-2){3}{\psline[arrows=->,arrowsize=3pt,linewidth=0.5pt](0,0)(0,-2.2)}\psline[linewidth=0.5pt](0,-7)(0,-8)}
\multirput(0,8)(2,0){5}{$\scriptstyle\bullet$}\multirput(0,6)(2,0){5}{$\scriptstyle\bullet$}\multirput(0,4)(2,0){5}{$\scriptstyle\bullet$}
\multirput(0,2)(2,0){5}{$\scriptstyle\bullet$}\multirput(0,0)(2,0){5}{$\scriptstyle\bullet$}
\rput[t](0,-0.3){$\scriptstyle0$}\rput[b](4.1,-0.9){$\cdots$}\rput[t](8.4,-0.3){$\scriptstyle n-1$}
\rput[r](-0.35,0){$\scriptstyle0$}\rput(-0.6,4.3){$\vdots$}\rput[r](-0.25,8){$\scriptstyle n-1$}
\multirput[b](1,8.15)(2,0){4}{$\scriptscriptstyle xz_1$}
\multirput[b](1,6.15)(2,0){4}{$\scriptscriptstyle xz_2$}
\multirput[b](1,4.25)(2,0){3}{$\scriptscriptstyle x$}\multirput[b](1,2.25)(2,0){2}{$\scriptscriptstyle x$}\rput[b](1,0.25){$\scriptscriptstyle x$}
\rput[b](7,4.2){$\scriptscriptstyle y$}\multirput[b](5,2.2)(2,0){2}{$\scriptscriptstyle y$}\multirput[b](3,0.2)(2,0){3}{$\scriptscriptstyle y$}
\multirput[l](0.2,7)(2,0){5}{$\scriptscriptstyle 1$}\multirput[l](0.2,5)(2,0){5}{$\scriptscriptstyle 1$}
\multirput[l](0.2,3)(2,0){5}{$\scriptscriptstyle 1$}\multirput[l](0.2,1)(2,0){5}{$\scriptscriptstyle 1$}
\endpspicture
\caption{The directed grid $\widetilde{\mathcal{G}}_n$, with associated edge weights.}\label{grid2}\end{figure}

Let $\widetilde{\mathcal{G}}_n$ be the $n\times n$ directed grid with
vertex set $\{(i,j)\mid i,j=0,\ldots,n-1\}$,
where $(i,j)$ is in the $i$th column from the left and $j$th row from the bottom, and for which there are
horizontal edges directed rightward from $(i,j)$ to $(i+1,j)$ for $i=0,\ldots,n-2$, $j=0,\ldots,n-1$, and
vertical edges directed downward from $(i,j)$ to $(i,j-1)$ for $i=0,\ldots,n-1$, $j=1,\ldots,n-1$.  This grid is shown
in Figure~\ref{grid2}, where the edge weights in this diagram will be introduced in Section~\ref{ZDPPsect}.
Note that, in contrast to $\mathcal{G}_n$ in Section~\ref{ASMbijsect}, the rows and columns of $\widetilde{\mathcal{G}}_n$
are labeled by Cartesian, rather than matrix-type, coordinates.

Now define
\begin{multline}\label{NILP1}\NILP(n)=\\
\left\{\begin{array}{@{}l@{}}\text{sets $P$ of noninter-}\\
\text{secting paths on }\widetilde{\mathcal{G}}_n\end{array}\,\left|\;\begin{array}{@{}l@{}}
\text{there exist }0\le t\le n-1\text{ and }n=\lambda_0>\lambda_1>\ldots>\\
\lambda_t>\lambda_{t+1}=0\text{ for which  }P\text{ consists of paths from }\\
(0,\lambda_{i-1}-1)\text{ to }(\lambda_i,0),\text{ for each }i=1,\ldots,t+1\end{array}\right.\right\}.\end{multline}
For example,
\psset{unit=5mm}
\begin{equation}\label{NILP3}\NILP(3)\,=\,\left\{\raisebox{-4mm}{
\pspicture(0,0)(2.8,2)\psgrid[subgriddiv=1,griddots=6,gridlabels=0](0,0)(2,2)
\psline[linewidth=1.7pt,linecolor=blue](0,2)(0,0)\rput(2.3,0.5){,}\endpspicture
\pspicture(0,0)(2.8,2)\psgrid[subgriddiv=1,griddots=6,gridlabels=0](0,0)(2,2)
\psline[linewidth=1.7pt,linecolor=blue](0,2)(2,2)(2,0)\psline[linewidth=1.7pt,linecolor=blue](0,1)(1,1)(1,0)
\psdots[dotscale=0.9,linecolor=blue](0,0)\rput(2.3,0.5){,}\endpspicture
\pspicture(0,0)(2.8,2)\psgrid[subgriddiv=1,griddots=6,gridlabels=0](0,0)(2,2)
\psline[linewidth=1.7pt,linecolor=blue](0,2)(0,1)(1,1)(1,0)
\psdots[dotscale=0.9,linecolor=blue](0,0)\rput(2.3,0.5){,}\endpspicture
\pspicture(0,0)(2.8,2)\psgrid[subgriddiv=1,griddots=6,gridlabels=0](0,0)(2,2)
\psline[linewidth=1.7pt,linecolor=blue](0,2)(2,2)(2,0)\psline[linewidth=1.7pt,linecolor=blue](0,1)(0,0)\rput(2.3,0.5){,}\endpspicture
\pspicture(0,0)(2.8,2)\psgrid[subgriddiv=1,griddots=6,gridlabels=0](0,0)(2,2)
\psline[linewidth=1.7pt,linecolor=blue](0,2)(1,2)(1,0)\psdots[dotscale=0.9,linecolor=blue](0,0)\rput(2.3,0.5){,}\endpspicture
\pspicture(0,0)(2.8,2)\psgrid[subgriddiv=1,griddots=6,gridlabels=0](0,0)(2,2)
\psline[linewidth=1.7pt,linecolor=blue](0,2)(1,2)(1,1)(2,1)(2,0)
\psline[linewidth=1.7pt,linecolor=blue](0,1)(0,0)\rput(2.3,0.5){,}\endpspicture
\pspicture(0,0)(2.2,2)\psgrid[subgriddiv=1,griddots=6,gridlabels=0](0,0)(2,2)
\psline[linewidth=1.7pt,linecolor=blue](0,2)(1,2)(1,0)(2,0)
\psline[linewidth=1.7pt,linecolor=blue](0,1)(0,0)\endpspicture}
\right\}.\end{equation}

Statistics for each $P\in\NILP(n)$ are now defined as
\begin{align}
\label{NILPstat1}\nu(P)&=\text{number of rightward steps of $P$ above the subdiagonal line }\{(i,i-1)\},\\
\mu(P)&=\text{number of rightward steps of $P$ below the subdiagonal line }\{(i,i-1)\},\\
\rho_1(P)&=\text{number of rightward steps of $P$ in row $n-1$ of }\widetilde{\mathcal{G}}_n,\\
\notag\rho_2(P)&=(\text{number of rightward steps of $P$ in row $n-2$ of }\widetilde{\mathcal{G}}_n)\;+\\
\label{NILPstat4}&\quad\;(\text{number of paths of $P$ which start at }(0,n-2)).\end{align}
Note that the second term on the RHS of \eqref{NILPstat4} (which is obviously 0 or 1) is also
the number of paths of $P$ which end at $(n-1,0)$.

It can be shown straightforwardly that there is a natural bijection between $\DPP(n)$ and $\NILP(n)$,
and that for each $D\in\ASM(n)$ and $P\in\NILP(n)$ which correspond under this bijection,
$\nu(D)=\nu(P)$, $\mu(D)=\mu(P)$, $\rho_1(D)=\rho_1(P)$ and $\rho_2(D)=\rho_2(P)$
(where the statistics in these equations are given
by~\eqref{nuD}--\eqref{rho2D} on the LHS and~\eqref{NILPstat1}--\eqref{NILPstat4} on the RHS).

The details of this bijection are as follows.
To map $D\in\DPP(n)$ to $P\in\NILP(n)$, first let~$t$ be the number of rows in~$D$, and let~$\lambda_i$ be the length
of row~$i$ of~$D$, as in~\eqref{DPP}.  Also define $\lambda_0=n$ and $\lambda_{t+1}=0$.  Then
obtain~$P$ by forming a path, for each $i=1,\ldots,t+1$, from $(0,\lambda_{i-1}-1)$ to
$(\lambda_i,0)$ whose~$k$th rightward step has height $D_{i,i+k-1}-1$.
To map $P\in\NILP(n)$ to $D\in\DPP(n)$,
set $D_{ij}$ to be $1$ plus the height of the $(j-i+1)$th rightward step in the~$i$th path from the top.

For example, the DPP~$D$ of~\eqref{ASMDPPEx} (with $n=6$) and the corresponding set of nonintersecting lattice paths are
\begin{equation}\label{NILPex}
\begin{array}{@{}c@{\:}c@{\:}c@{\:}c@{\:}c@{\:}c@{}}
&6&6&6&5&2\\
&&4&4&1\\
&&&3\end{array}\quad\longleftrightarrow\quad
\psset{unit=7.5mm}\raisebox{-20mm}{\pspicture(-0.7,-0.4)(5,5.3)
\psgrid[subgriddiv=1,griddots=9,gridlabels=0](0,0)(5,5)
\psline[linewidth=1pt](0.9,-0.1)(5.1,4.1)
\psline[linewidth=2.5pt,linecolor=blue](0,5)(3,5)(3,4)(4,4)(4,1)(5,1)(5,0)
\psline[linewidth=2.5pt,linecolor=blue](0,4)(0,3)(2,3)(2,0)(3,0)
\psline[linewidth=2.5pt,linecolor=blue](0,2)(1,2)(1,0)
\psdots[dotscale=1.5,linecolor=blue](0,0)
\rput[b](0.5,5.15){$\color{red}\scriptstyle6$}\rput[b](1.5,5.15){$\color{red}\scriptstyle6$}
\rput[b](2.5,5.15){$\color{red}\scriptstyle6$}\rput[b](3.5,4.15){$\color{red}\scriptstyle5$}
\rput[b](4.5,1.15){$\color{green}\scriptstyle2$}
\rput[b](0.5,3.15){$\color{red}\scriptstyle4$}\rput[b](1.5,3.15){$\color{red}\scriptstyle4$}
\rput[b](0.5,2.15){$\color{red}\scriptstyle3$}\rput[b](2.5,0.15){$\color{green}\scriptstyle1$}
\rput[t](0,-0.2){$\scriptscriptstyle0$}\rput[t](1,-0.2){$\scriptscriptstyle1$}\rput[t](2,-0.2){$\scriptscriptstyle2$}
\rput[t](3,-0.2){$\scriptscriptstyle3$}\rput[t](4,-0.2){$\scriptscriptstyle4$}\rput[t](5,-0.2){$\scriptscriptstyle5$}
\rput[r](-0.2,0){$\scriptscriptstyle0$}\rput[r](-0.2,1){$\scriptscriptstyle1$}\rput[r](-0.2,2){$\scriptscriptstyle2$}
\rput[r](-0.2,3){$\scriptscriptstyle3$}\rput[r](-0.2,4){$\scriptscriptstyle4$}\rput[r](-0.2,5){$\scriptscriptstyle5$}
\endpspicture}\end{equation}
where the subdiagonal line $\{(i,i-1)\}$ is shown,
and each part of $D$ is displayed above its corresponding rightward step, with the
nonspecial and special parts in red and green respectively.  It can be seen that the statistics~\eqref{NILPstat1}--\eqref{NILPstat4} for the
set of paths $P$ in~\eqref{NILPex} are $\nu(P)=7$, $\mu(P)=2$, $\rho_1(P)=3$ and $\rho_2(P)=2$, and that these match the
respective statistics~\eqref{nuD}--\eqref{rho2D} for the DPP~$D$ in~\eqref{NILPex}, as given in~\eqref{ASMDPPstatex}.

As further examples, in~\eqref{ASMDPP3} and~\eqref{NILP3} the elements of $\DPP(3)$ and $\NILP(3)$
are listed in an order for which respective elements correspond under the bijection of this section.
It can be seen that $Z^\DPP_3(x,y,z_1,z_2)$, as given in~\eqref{ZASMDPP3},
could now be obtained using this bijection and~\eqref{NILP3}.

\subsection{The Lindstr\"{o}m--Gessel--Viennot theorem}\label{LGVsect}
In this section, the Lindstr\"{o}m--Gessel--Viennot theorem for the weighted enumeration of sets of nonintersecting
paths in terms of a determinant is stated.

Consider an acyclic directed graph~$G$, let a weight be assigned to each edge of~$G$, and
define the weight~$W(p)$ of a path~$p$ on~$G$ to be the product of the weights of
the edges along which~$p$ passes.
For vertices~$u$ and~$v$ of~$G$, let $\mathcal{P}_{u,v}$ denote the set of all
paths on~$G$ from~$u$ to~$v$.   For vertices $u_1,\ldots,u_m$ and $v_1,\ldots,v_m$ of~$G$,
let $\mathcal{N}_{G}(u_1,\ldots,u_m;v_1,\ldots,v_m)$ denote the set of all sets~$P$ of
paths on~$G$ such that $P$ consists of a path of $\mathcal{P}_{u_i,v_i}$ for each $i=1,\ldots,m$, and
different paths of~$P$ are nonintersecting.
The Lindstr\"{o}m--Gessel--Viennot theorem~\cite{GesVie85,GesVie89,Lin73} can now be stated as follows.
\begin{theorem*}[Lindstr\"{o}m; Gessel, Viennot]
If $\mathcal{N}_{G}(u_{\sigma_1},\ldots,u_{\sigma_m};v_1,\ldots,v_m)$ is empty for each permutation $\sigma$
of $\{1,\ldots,m\}$ other than the identity, then
\begin{equation}\label{LGV}\sum_{P\in\mathcal{N}_{G}(u_1,\ldots,u_m;v_1,\ldots,v_m)}\,
\prod_{p\in P}W(p)\,=\det_{1\le i,j\le m}\Biggl(\,
\sum_{p\in\mathcal{P}_{u_i,v_j}}\!W(p)\Biggr).\end{equation}
\end{theorem*}
For a proof of~\eqref{LGV} see, for example, Gessel and Viennot~\cite[Sec.~2]{GesVie89}
or Stembridge~\cite[Sec.~1]{Ste90}.

\subsection{Determinant formula for the DPP doubly-refined generating function}\label{ZDPPsect}
In this section, the bijection of Section~\ref{DPPbijsect} and the Lindstr\"{o}m--Gessel--Viennot theorem~\eqref{LGV} are
used to obtain a determinant formula, as already stated in the second equality of~\eqref{doubdet},
 for the doubly-refined DPP generating function~\eqref{ZDPP}.

The directed graph $G$ of Section~\ref{LGVsect} is now taken to be the $n\times n$ directed grid $\widetilde{\mathcal{G}}_n$,
defined in Section~\ref{DPPbijsect}.

The set $\NILP(n)$, as defined in~\eqref{NILP1}, can be written using the
notation of Section~\ref{LGVsect} as
\begin{multline}\label{NILP2}
\NILP(n)=\\
\bigcup_{\substack{0\le t\le n-1\\\rule{0ex}{1.5ex}
n-1\ge\lambda_1>\ldots>\lambda_t\ge1}}\!
\mathcal{N}_{\widetilde{\mathcal{G}}_n}\bigl((0,n-1),(0,\lambda_1-1),\ldots,(0,\lambda_t-1);
(\lambda_1,0),\ldots,(\lambda_t,0),(0,0)\bigr).\end{multline}

Now assign weights to the edges of $\widetilde{\mathcal{G}}_n$, as indicated in Figure~\ref{grid2}.
More specifically, assign~$xz_1$ to each horizontal edge
in row $n-1$,~$xz_2$ to each horizontal edge in
row~$n-2$,~$x$ to each horizontal edge from $(i,j)$ to $(i+1,j)$ for $0\le i\le j\le n-3$,~$y$ to each horizontal edge
from $(i,j)$ to $(i+1,j)$ for $0\le j<i\le n-2$, and~1 to each vertical edge.
Also, slightly expanding the notation of Section~\ref{LGVsect},
denote the weight of a path $p$ on $\widetilde{\mathcal{G}}_n$ as $W(x,y,z_1,z_2,p)$.

The bijection of Section~\ref{DPPbijsect} and the Lindstr\"{o}m--Gessel--Viennot theorem~\eqref{LGV}, together
with~\eqref{ZDPP},~\eqref{NILPstat1}--\eqref{NILPstat4},~\eqref{NILP2}
and the observation that the condition for the validity of~\eqref{LGV} is satisfied,
now give
\begin{multline}\label{DPPLGV1}
Z^\DPP_n(x,y,z_1,z_2)=\\[-2mm]
\sum_{\substack{0\le t\le n-1\\\rule{0ex}{1.5ex}
1\le\lambda_t<\ldots<\lambda_1\le n-1}}\left(z_2^{\;\delta_{\lambda_1-1,n-2}}\!\!
\det_{\rule{0ex}{2.3ex}\substack{i=0,\lambda_t,\ldots,\lambda_1\\j=\lambda_t-1,\ldots,\lambda_1-1,n-1}}
\Biggl(\,\sum_{p\in\mathcal{P}_{(0,j),(i,0)}}W(x,y,z_1,z_2,p)\Biggr)\right).
\end{multline}
Note that the term $z_2^{\;\delta_{\lambda_1-1,n-2}}$ ($=z_2^{\;\delta_{\lambda_1,n-1}}$) on the
RHS of~\eqref{DPPLGV1} arises from the second term on the RHS of~\eqref{NILPstat4}.

It can be shown straightforwardly that, for any matrix $\bigl(M_{ij}\bigr)_{0\le i,j\le n-1}$,
\begin{equation}\label{detidentity}
\det_{0\le i,j\le n-1}\bigl(M_{ij}-\delta_{i,j+1}\bigr)=
\sum_{S\subset\{1,\ldots,n-1\}}\det M_{\{0\}\cup S,(S-1)\cup\{n-1\}},\end{equation}
where  $M_{\{0\}\cup S,(S-1)\cup\{n-1\}}$ denotes the
submatrix of $M$ formed by restricting the rows and columns to those indexed by
$\{0\}\cup S$ and $\{s-1\mid s\in S\}\cup\{n-1\}$ respectively.

Applying~\eqref{detidentity} to~\eqref{DPPLGV1},
and taking account of the term $z_2^{\;\delta_{\lambda_1-1,n-2}}$, gives
\begin{equation}\label{DPPLGV}
Z^\DPP_n(x,y,z_1,z_2)=
\det_{0\le i,j\le n-1}\Biggl(-\delta_{i,j+1}+
z_2^{\;\delta_{j,n-2}}\!\!\sum_{p\in\mathcal{P}_{(0,j),(i,0)}}W(x,y,z_1,z_2,p)\Biggr).\end{equation}
(Note that $Z^\DPP_n(x,y,z_1,z_2)=
\det_{0\le i,j\le n-1}\bigl(-\delta_{i,j+1}+
z_2^{\;\delta_{i,n-1}}$ $\sum_{p\in\mathcal{P}_{(0,j),(i,0)}}\!W(x,y,z_1,z_2,p)\bigr)$ is also valid.)

It is now found that
\begin{equation}\label{DPPwp}
\sum_{p\in\mathcal{P}_{(0,j),(i,0)}}\!\!\!W(x,y,z_1,z_2,p)=
\begin{cases}
\sum_{k=0}^{\min(i,j+1)}\binom{i-1}{i-k}\binom{j+1}{k}x^ky^{i-k},&j\le n-3,\\
\sum_{k=0}^i\sum_{l=0}^k\binom{i-1}{i-k}\binom{n-l-2}{k-l}x^ky^{i-k}z_2^l,&j=n-2,\\
\sum_{k=0}^i\sum_{l=0}^k\sum_{m=0}^l\binom{i-1}{i-k}\binom{n-l-2}{k-l}x^ky^{i-k}z_1^mz_2^{l-m},&j=n-1.
\end{cases}
\end{equation}
More specifically, the sums of weights of paths of $\mathcal{P}_{(0,j),(i,0)}$ in~\eqref{DPPwp} can be obtained as follows,
with these derivations being shown diagrammatically in Figure~\ref{DPPdetfig}.
For $j\le n-3$, combine any
of the $\binom{j+1}{k}$ paths of $\mathcal{P}_{(0,j),(k,k-1)}$, each with weight~$x^k$,
and any of the $\binom{i-1}{i-k}$ paths of $\mathcal{P}_{(k,k-1),(i,0)}$, each with weight~$y^{i-k}$,
for any $0\le k\le\min(i,j+1)$.
For $j=n-2$, combine the single path of $\mathcal{P}_{(0,n-2),(l,n-2)}$, with
weight~$(xz_2)^l$, the single path of $\mathcal{P}_{(l,n-2),(l,n-3)}$,
with weight~$1$, any of the $\binom{n-l-2}{k-l}$ paths of
$\mathcal{P}_{(l,n-3),(k,k-1)}$, each with weight~$x^{k-l}$, and
any of the $\binom{i-1}{i-k}$ paths of $\mathcal{P}_{(k,k-1),(i,0)}$, each
with weight~$y^{i-k}$, for any $0\le l\le k\le i$.
For $j=n-1$, combine the single path of $\mathcal{P}_{(0,n-1),(m,n-1)}$, with
weight~$(xz_1)^m$, the single path of $\mathcal{P}_{(m,n-1),(m,n-2)}$, with
weight~$1$, the single path of $\mathcal{P}_{(m,n-2),(l,n-2)}$, with
weight~$(xz_2)^{l-m}$, the single path of $\mathcal{P}_{(l,n-2),(l,n-3)}$,
with weight~$1$, any of the $\binom{n-l-2}{k-l}$ paths of
$\mathcal{P}_{(l,n-3),(k,k-1)}$, each with weight~$x^{k-l}$, and
any of the $\binom{i-1}{i-k}$ paths of $\mathcal{P}_{(k,k-1),(i,0)}$, each
with weight~$y^{i-k}$, for any $0\le m\le l\le k\le i$.

\begin{figure}[h]\centering
\psset{unit=4.2mm}
\pspicture(-1.7,-1.8)(10.7,10)
\psgrid[subgriddiv=1,griddots=6,gridlabels=0](0,0)(10,10)
\psline[linewidth=1pt](0.9,-0.1)(10.1,9.1)
\rput(0,-0.55){$\scriptstyle0$}\rput(4,-0.55){$\scriptstyle k$}\rput(8,-0.55){$\scriptstyle i$}\rput(10.1,-0.55){$\scriptstyle n-1$}
\rput[r](-0.3,0){$\scriptstyle0$}\rput[r](-0.2,3){$\scriptstyle k-1$}\rput[r](-0.4,6){$\scriptstyle j$}
\rput[r](-0.28,9){$\scriptstyle n-2$}\rput[r](-0.28,8){$\scriptstyle n-3$}\rput[r](-0.2,10){$\scriptstyle n-1$}
\rput(5,-2){$j\le n-3$}
\psline[linewidth=1.7pt,linecolor=blue](0,6)(3,6)(3,3)(4,3)(4,2)(6,2)(6,1)(7,1)(7,0)(8,0)
\psdots[dotscale=1.2](0,6)(4,3)(8,0)
\rput[b](0.5,6.2){$\scriptscriptstyle x$}\rput[b](1.5,6.2){$\scriptscriptstyle x$}
\rput[b](2.5,6.2){$\scriptscriptstyle x$}\rput[b](3.5,3.2){$\scriptscriptstyle x$}\rput[b](4.5,2.2){$\scriptscriptstyle y$}
\rput[b](5.5,2.2){$\scriptscriptstyle y$}\rput[b](6.5,1.2){$\scriptscriptstyle y$}\rput[b](7.5,0.2){$\scriptscriptstyle y$}
\endpspicture\quad
\pspicture(-1.7,-1.8)(10.7,10)
\psgrid[subgriddiv=1,griddots=6,gridlabels=0](0,0)(10,10)
\psline[linewidth=1pt](0.9,-0.1)(10.1,9.1)
\rput(0,-0.55){$\scriptstyle0$}\rput(5,-0.55){$\scriptstyle k$}\rput(2,-0.55){$\scriptstyle l$}
\rput(8,-0.55){$\scriptstyle i$}\rput(10.1,-0.55){$\scriptstyle n-1$}
\rput[r](-0.3,0){$\scriptstyle0$}\rput[r](-0.2,4){$\scriptstyle k-1$}
\rput[r](-0.28,9){$\scriptstyle n-2$}\rput[r](-0.28,8){$\scriptstyle n-3$}\rput[r](-0.2,10){$\scriptstyle n-1$}
\rput(5,-2){$j=n-2$}
\psline[linewidth=1.7pt,linecolor=blue](0,9)(2,9)(2,6)(4,6)(4,5)(5,5)(5,3)(6,3)(6,2)(8,2)(8,0)
\psdots[dotscale=1.2](0,9)(2,9)(2,8)(5,4)(8,0)
\rput[b](0.51,9.23){$\scriptscriptstyle x\!z_2$}\rput[b](1.52,9.23){$\scriptscriptstyle x\!z_2$}
\rput[b](2.5,6.2){$\scriptscriptstyle x$}\rput[b](3.5,6.2){$\scriptscriptstyle x$}\rput[b](4.5,5.2){$\scriptscriptstyle x$}
\rput[b](5.5,3.2){$\scriptscriptstyle y$}\rput[b](6.5,2.2){$\scriptscriptstyle y$}
\rput[b](7.5,2.2){$\scriptscriptstyle y$}
\endpspicture\quad
\pspicture(-1.7,-1.8)(10.7,10)
\psgrid[subgriddiv=1,griddots=6,gridlabels=0](0,0)(10,10)
\psline[linewidth=1pt](0.9,-0.1)(10.1,9.1)
\rput(0,-0.55){$\scriptstyle0$}\rput(3,-0.55){$\scriptstyle m$}\rput(5,-0.55){$\scriptstyle l$}
\rput(7,-0.55){$\scriptstyle k$}\rput(8,-0.55){$\scriptstyle i$}\rput(10.1,-0.55){$\scriptstyle n-1$}
\rput[r](-0.3,0){$\scriptstyle0$}\rput[r](-0.2,6){$\scriptstyle k-1$}
\rput[r](-0.28,9){$\scriptstyle n-2$}\rput[r](-0.28,8){$\scriptstyle n-3$}\rput[r](-0.2,10){$\scriptstyle n-1$}
\rput(5,-2){$j=n-1$}
\psline[linewidth=1.7pt,linecolor=blue](0,10)(3,10)(3,9)(5,9)(5,6)(7,6)(7,1)(8,1)(8,0)
\psdots[dotscale=1.2](0,10)(3,10)(3,9)(5,9)(5,8)(7,6)(8,0)
\rput[b](0.45,10.2){$\scriptscriptstyle xz_1$}\rput[b](1.5,10.22){$\scriptscriptstyle xz_1$}\rput[b](2.56,10.22){$\scriptscriptstyle xz_1$}
\rput[b](3.51,9.23){$\scriptscriptstyle x\!z_2$}\rput[b](4.51,9.23){$\scriptscriptstyle x\!z_2$}
\rput[b](5.5,6.2){$\scriptscriptstyle x$}\rput[b](6.5,6.2){$\scriptscriptstyle x$}
\rput[b](7.5,1.2){$\scriptscriptstyle y$}
\endpspicture
\caption{Derivation of~\eqref{DPPwp}.}\label{DPPdetfig}\end{figure}
Finally, using~\eqref{DPPwp} in~\eqref{DPPLGV} gives the second equality of~\eqref{doubdet}, i.e.,
\begin{equation}\label{ZDPPdet1}Z^\DPP_n(x,y,z_1,z_2)=\det_{0\le i,j\le n-1}\bigl(K_n(x,y,z_1,z_2)_{ij}\bigr),\end{equation}
where $K_n(x,y,z_1,z_2)_{ij}$ is defined in~\eqref{K}.

\subsection{Transformation of the DPP determinant formula}
In this section, a further determinant formula for the DPP doubly-refined generating function~\eqref{ZDPP}
is obtained by elementary transformation of the determinant formula~\eqref{ZDPPdet1}.

Define, for $i=0,\ldots,n-1$,
\begin{align}
\notag C_n(x,y,z)_i&=\textstyle\sum_{k=0}^i\sum_{l=0}^k\binom{i-1}{i-k}\binom{n-l-2}{k-l}x^ky^{i-k}z^{l+1},\\
\label{C}C_n(x,y,z_1,z_2)_i&=\textstyle \sum_{k=0}^i\sum_{l=0}^k\sum_{m=0}^l\binom{i-1}{i-k}
\binom{n-l-2}{k-l}x^ky^{i-k}z_1^mz_2^{l-m}.
\end{align}
Thus, $C_n(x,y,z_2)_i-\delta_{i,n-1}$ and $C_n(x,y,z_1,z_2)_i$ are the entries in row~$i$ of the
second-last and last columns respectively of the matrix $K_n(x,y,z_1,z_2)$, as
defined in~\eqref{K} and used in~\eqref{ZDPPdet1}.

It is now found that
\begin{equation}\label{Ceq}(z_1-z_2)C_n(x,y,z_1,z_2)_i=C_n(x,y,z_1)_i-C_n(x,y,z_2)_i.\end{equation}
This can be proved combinatorially as follows. As shown in Section~\ref{ZDPPsect},
and stated in the cases $j=n-1$ and $j=n-2$ of~\eqref{DPPwp},
$C_n(x,y,z_1,z_2)_i=\sum_{p\in\mathcal{P}_{(0,n-1),(i,0)}}\!\!W(x,y,z_1,z_2,p)$ and
$z_2{}^{\!-1}\,C_n(x,y,z_2)_i=\sum_{p\in\mathcal{P}_{(0,n-2),(i,0)}}\!\!W(x,y,z_1,z_2,p)$.
By partitioning $\mathcal{P}_{(0,n-1),(i,0)}$ into those paths with an initial downward step,
which correspond (by deleting the initial step, with weight 1) to $\mathcal{P}_{(0,n-2),(i,0)}$,
and those paths with an initial rightward step, which correspond (by deleting the initial step, with weight $xz_1$)
to $\mathcal{P}_{(1,n-1),(i,0)}$, it follows that
\begin{equation}\label{Ceq1}(xz_1){}^{-1}\bigl(C_n(x,y,z_1,z_2)_i-z_2{}^{\!-1}\,C_n(x,y,z_2)_i\bigr)=
\sum_{p\in\mathcal{P}_{(1,n-1),(i,0)}}W(x,y,z_1,z_2,p).\end{equation}
Now consider the involutions on $\mathcal{P}_{(0,n-1),(i,0)}$ and on $\mathcal{P}_{(1,n-1),(i,0)}$ in which
a path which starts with $k_1$ rightward steps in row $n-1$ of $\widetilde{\mathcal{G}}_n$,
followed by $k_2$ rightward steps in row~$n-2$ of~$\widetilde{\mathcal{G}}_n$, is
mapped to the path in which the roles of $k_1$ and $k_2$ are interchanged, with the part of the path below row $n-2$ being
kept the same. Using these involutions, it follows that both $C_n(x,y,z_1,z_2)$ and
the RHS of~\eqref{Ceq1} are symmetric in $z_1$ and $z_2$.  Finally, interchanging $z_1$ and $z_2$ on the LHS of~\eqref{Ceq1},
and setting this equal to the original LHS of~\eqref{Ceq1}, gives~\eqref{Ceq}.

Alternatively,~\eqref{Ceq} can be proved algebraically as follows.
Let $[u^i]f(u)$ denote the coefficient of~$u^i$ in a power series~$f(u)$.
Then, briefly, $C_n(x,y,z)_i=[u^i]\,\frac{z(1-yu)}{1-(x(z-1)+y)u}\bigl(1+\frac{xu}{1-yu}\bigr)^{n-1}$, and
$C_n(x,y,z_1)_i-C_n(x,y,z_2)_i=[u^i]\,\bigl(\frac{z_1}{1-(x(z_1-1)+y)u}-
\frac{z_2}{1-(x(z_2-1)+y)u}\bigr)(1-yu)\bigl(1+\frac{xu}{1-yu}\bigr)^{n-1}
=[u^i]\,\frac{(z_1-z_2)(1-yu)^2}{(1-(x(z_1-1)+y)u)(1-(x(z_2-1)+y)u)}\bigl(1+\frac{xu}{1-yu}\bigr)^n
=(z_1-z_2)C_n(x,y,z_1,z_2)_i$.

Starting with the determinant formula~\eqref{ZDPPdet1}, together with definitions~\eqref{K} and~\eqref{C},
and then successively multiplying the last column
of $K_n(x,y,z_1,z_2)$
by $z_1-z_2$, adding the second-last column to the last column,
applying~\eqref{Ceq}, and interchanging the last two columns,
it now follows that
\begin{equation}\label{Ldet}(z_2-z_1)Z^\DPP_n(x,y,z_1,z_2)=\det_{0\le i,j\le n-1}\bigl(L_n(x,y,z_1,z_2)_{ij}\bigr),\end{equation}
where
\begin{equation}\label{L}L_n(x,y,z_1,z_2)_{ij}=
\begin{cases}
\sum_{k=0}^{\min(i,j+1)}\binom{i-1}{i-k}\binom{j+1}{k}x^ky^{i-k}-\delta_{i,j+1},&j\le n-3,\\
C_n(x,y,z_1)_i-\delta_{i,n-1},&j=n-2,\\
C_n(x,y,z_2)_i-\delta_{i,n-1},&j=n-1.\\
\end{cases}
\end{equation}

\subsection{Application of the Desnanot--Jacobi identity}
In this section, the final steps in the proof of the DPP case of~\eqref{propeq2} are taken.

Define an $n\times(n+2)$ matrix $\bigl(L_n(x,y,z_1,z_2,z_3,z_4)_{ij}\bigr)_{0\le i\le n-1;\,0\le j\le n+1}$ by
\begin{equation}\label{LL}L_n(x,y,z_1,z_2,z_3,z_4)_{ij}=\begin{cases}L_n(x,y,z_1,z_2)_{ij},&j\le n-1,\\
C_n(x,y,z_3)_i-\delta_{i,n-1},&j=n,\\
C_n(x,y,z_4)_i-\delta_{i,n-1},&j=n+1,\end{cases}\end{equation}
where $L_n(x,y,z_1,z_2)_{ij}$ and $C_n(x,y,z)_i$ are given by~\eqref{L} and~\eqref{C}.

Applying the form~\eqref{Bazin} of the Desnanot--Jacobi identity with its four selected rows taken as
$k_i=n+i-3$ (and with row and column numbers now starting from~$0$) to the
transpose of the matrix~\eqref{LL},
and using the determinant formula~\eqref{Ldet} on each of the six determinants which appear,
now gives the DPP case of~\eqref{propeq2}.

\section{Discussion}\label{discuss}
The proof of Theorem~\ref{DRT}, which has been the primary focus of this paper, has now been presented.
In summary, the proof involved
using certain bijections together with the Izergin--Korepin formula~\eqref{Idet1} and Lindstr\"{o}m--Gessel--Viennot theorem~\eqref{LGV} to obtain
determinantal expressions for the doubly-refined enumeration of ASMs and DPPs,
then applying the form~\eqref{Bazin} of the Desnanot--Jacobi identity
to obtain the identity~\eqref{propeq2} satisfied by the ASM and DPP doubly-refined generating functions,
then using the argument given in Section~\ref{mainres} (i.e., setting~$z_3=1$ and $z_4=0$ in~\eqref{propeq2} and applying~\eqref{Zsing}--\eqref{Z0})
to show that the doubly-refined generating functions can be expressed in terms of their singly-refined counterparts by~\eqref{propeq1},
and finally using the equality~\eqref{MRR} of the singly-refined generating functions to establish the required
equality~\eqref{DRT1} of the doubly-refined generating functions.

A variation of this proof involves applying the form~\eqref{DesJacBazin}, instead of~\eqref{Bazin},
of the Desnanot--Jacobi identity to the ASM and DPP determinants,
which leads directly to~\eqref{propeq1}, instead of~\eqref{propeq2}.  However,
in this approach, different terms of~\eqref{propeq1} need to be treated separately at certain stages of the derivation,
whereas the approach which was used seems preferable, since all the terms of~\eqref{propeq2}
can be treated uniformly throughout.

In the remainder of this section, some further matters related to the main results of this paper are discussed.

\subsection{Special cases}
In this section, some special cases of Theorem~\ref{DRT} which can be proved bijectively
are discussed briefly.

It follows from Theorem~\ref{DRT} that, for each~$n$, there exist bijections between $\ASM(n)$ and $\DPP(n)$ with the property that,
for each $A\in\ASM(n)$ and $D\in\DPP(n)$ which correspond under the bijection,
$\nu(A)=\nu(D)$, $\mu(A)=\mu(D)$, $\rho_1(A)=\rho_1(D)$ and $\rho_2(A)=\rho_2(D)$.

However, although these bijections necessarily exist, a natural and explicit such bijection for arbitrary~$n$ is not currently known.
Nevertheless, explicit bijections are known
for some special classes of the subsets of $\ASM(n)$ and $\DPP(n)$ in Theorem~\ref{DRT}
(i.e., for certain classes of values of~$p$,~$m$,~$k_1$ and~$k_2$ in Theorem~\ref{DRT}), with these bijections
thereby providing alternative proofs of
certain cases of the theorem. Some of these cases will now be outlined, although most of the details will be omitted.

For $k_1=k_2$, $p=k_1(k_1+1)/2$ and $m=k_1(n-k_1-1)$,
the sets in Theorem~\ref{DRT} consist of
a single ASM and single DPP, for each~$0\le k_1=k_2\le n-1$ and arbitrary~$n$, so that there is a single, trivial bijection.
The ASM and DPP are given explicitly in~\cite[Eqs.~(21)--(22)]{BehDifZin12a}.

For $m=0$, with $n$, $p$, $k_1$ and $k_2$ arbitrary (i.e., for permutation matrices and DPPs with no special parts),
a bijection can be obtained by modifying a known bijection
for the corresponding singly-refined case in which~$\rho_2$ and~$k_2$ do not appear,
where the modification is needed to give the property that $\rho_2(A)=\rho_2(D)$ ($=k_2$) for any
corresponding~ASM $A$ and DPP~$D$.
Details and references for the known bijection are given in~\cite[pp.~336--337]{BehDifZin12a}.
It can also be shown straightforwardly for this case that
\begin{multline}Z^\ASM_n(x,0,z_1,z_2)=Z^\DPP_n(x,0,z_1,z_2)=\\
[n-2]_x!\sum_{0\le i<j\le n-1}(x^{n+i-j-1}z_1^iz_2^{n-j-1}+x^{n-i+j-2}z_1^{n-i-1}z_2^j),\end{multline}
where $[n]_x=1+x+\ldots+x^{n-1}$ and $[n]_x!=[n]_x[n-1]_x\ldots[1]_x$.

For $p=1$, with $n$, $m$, $k_1$ and $k_2$ arbitrary,
a bijection can again be obtained by modifying a known bijection
for the corresponding singly-refined case in which~$\rho_2$ and~$k_2$ do not appear,
with the modification needed to ensure that $\rho_2(A)=\rho_2(D)$ for any
corresponding~ASM $A$ and DPP~$D$.
Details for the known bijection are essentially given in~\cite[p.~337]{BehDifZin12a}.
This case, together with the operation~$\ast$ which will be described in Section~\ref{symm}, can also be used to give a bijection
for the case $p=n(n-1)/2-m-1$, with $n$,~$m$,~$k_1$ and~$k_2$ arbitrary.

In these cases, the bijections for the singly-refined cases seem reasonably natural, whereas the subsequent modifications
do not. In fact, it seems that
a natural, full bijection between $\ASM(n)$ and $\DPP(n)$ might satisfy
$\nu(A)=\nu(D)$, $\mu(A)=\mu(D)$, $\rho_1(A)=\rho_1(D)$ and $\rho_2(A)=\rho_2(\sigma(D))$ for each corresponding $A\in\ASM(n)$
and $D\in\DPP(n)$, where~$\sigma$ is a permutation (other than the identity) of $\DPP(n)$ which satisfies $\nu(\sigma(D))=\nu(D)$,
$\mu(\sigma(D))=\mu(D)$ and $\rho_1(\sigma(D))=\rho_1(D)$ for each $D\in\DPP(n)$.  If this is so, then it would be natural to replace
the fourth statistic for each $D\in\DPP(n)$ by $\rho_2(\sigma(D))$.

\subsection{Symmetry operations}\label{symm}
In this section, two symmetry operations on $\ASM(n)$ and $\DPP(n)$ which lead directly to some common properties
of the doubly-refined generating functions $Z^\ASM_n(x,y,z_1,z_2)$ and $Z^\DPP_n(x,y,z_1,z_2)$ are discussed.

Operations $\ast$ and $\dagger$ which transform $A\in\ASM(n)$ to $A^\ast,A^\dagger\in\ASM(n)$ are defined by
\begin{equation}A^*_{ij}=A_{i,n+1-j},\qquad A^\dagger_{ij}=A_{n+1-i,n+1-j}.\end{equation}
Thus, $\ast$ and $\dagger$ correspond to reflection of $A$ in a vertical line,
and rotation of~$A$ by~$\pi$ respectively.  These operations, together with the other natural reflections and rotations
of ASMs, were first considered by Mills, Robbins and
Rumsey~\cite{MilRobRum83},~\cite{MilRobRum86},~\cite{Rob91},~\cite{Rob00}, who formulated numerous conjectures involving the sizes of
sets of ASMs invariant under groups of such operations.  Many of these conjectures have since been proved,
and further related results have been obtained, by, for example,
Kuperberg~\cite{Kup02}, Okada~\cite{Oka06}, and Razumov and Stroganov~\cite{RazStr06b,RazStr06a}.

An operation $\ast$ which transforms $D\in\DPP(n)$ to $D^\ast\in\DPP(n)$ is defined as follows.
If $D_{n-j,n-i}$ is defined and $D_{n-j,n-i}\le j-i$ then $D^*_{ij}=j-i+1-D_{n-j,n-i}$, if
$i\le j<n$ and $D_{n-j,n-i}$ is not defined then $D^*_{ij}=n+1-i-($number of positive integers~$k$
satisfying $n+2-i-D_{k,n-j}\le k\le n-j)$, and if $D_{n-j,n-i}>j-i$ then $D^*_{ij}$ is not defined.
This operation was first defined by Mills, Robbins and Rumsey~\cite[p.~351]{MilRobRum83},
using the previous rather complicated description.
However, it was shown by Lalonde~\cite{Lal03} that the operation has a much simpler description
in terms of the set of nonintersecting lattice paths which correspond to a DPP.
(Lalonde used slightly different nonintersecting paths from the paths of~\eqref{NILP1} used here, but it can be shown
straightforwardly that~$\ast$ also has a relatively simple description in terms of~\eqref{NILP1}.)
Furthermore, it was shown by Krattenthaler~\cite{Kra06} that there is a natural bijection between $\DPP(n)$ and
the set of cyclically symmetric
tilings with three types of unit rhombi
of a certain punctured hexagon with alternating sides of length $n-1$ and~$n+1$,
and that, in terms of such tilings, $\ast$ simply corresponds to the reflection of a tiling in a symmetry axis of the hexagon.

An operation $\dagger$ which transforms $D\in\DPP(n)$ to $D^\dagger\in\DPP(n)$ is defined as follows.
Let~$\rho_1(D)$ and~$\rho_2(D)$ be given by~\eqref{rho1D}--\eqref{rho2D} and~$\rho_3(D)$ be the number
of $(n-1)$'s in the first row of $D$.  Then obtain $D^\dagger$ by replacing the~$\rho_1(D)$~$n$'s and $\rho_3(D)$ $(n-1)$'s
in the first row of~$D$ by $\rho_2(D)$ $n$'s followed by $(\rho_1(D)+\rho_3(D)-\rho_2(D))$ $(n-1)$'s, leaving all other parts of~$D$ unchanged.
This operation has not previously been considered in the literature.

It can be seen that the operations $\ast$ and $\dagger$ are involutions on $\ASM(n)$ and on $\DPP(n)$.
Furthermore, the statistics~\eqref{muA}--\eqref{rho2D} behave under these operations according to
\begin{align}\notag\nu(X^\ast)&=\tfrac{n(n-1)}{2}-\nu(X)-\mu(X),&\mu(X^\ast)&=\mu(X),\\[2mm]
\label{aststat}\rho_1(X^\ast)&=n-1-\rho_1(X),&\rho_2(X^\ast)&=n-1-\rho_2(X),\end{align}
and
\begin{equation}\label{dagstat}\nu(X^\dagger)=\nu(X),\qquad\mu(X^\dagger)=\mu(X),\qquad
\rho_1(X^\dagger)=\rho_2(X),\end{equation}
for any $X\in\ASM(n)$ or $X\in\DPP(n)$.  Proofs of the first three cases of~\eqref{aststat} for ASMs and DPPs
are given by Mills, Robbins and Rumsey~\cite[p.~352]{MilRobRum83}, and for DPPs
by Krattenthaler~\cite[pp.~1143--1144]{Kra06} and Lalonde~\cite[p.~317]{Lal03}.
The remaining cases in~\eqref{aststat}--\eqref{dagstat} can be obtained similarly.

It follows from~\eqref{aststat} that~$\ast$, when restricted to the set of ASMs or the set of DPPs in Theorem~\ref{DRT}, provides
a bijection to the corresponding set in which~$p$, $k_1$ and $k_2$ are replaced by~$n(n-1)/2-p-m$,~$n-1-k_1$ and~$n-1-k_2$ respectively.
Therefore,
\begin{align}\notag Z^\ASM_n(x,y,z_1,z_2)&=x^{n(n-1)/2}\,(z_1z_2)^{n-1}\,Z^\ASM_n(\tfrac{1}{x},\tfrac{y}{x},\tfrac{1}{z_1},\tfrac{1}{z_2}),\\
Z^\DPP_n(x,y,z_1,z_2)&=x^{n(n-1)/2}\,(z_1z_2)^{n-1}\,Z^\DPP_n(\tfrac{1}{x},\tfrac{y}{x},\tfrac{1}{z_1},\tfrac{1}{z_2}).\end{align}

Similarly, it follows from~\eqref{dagstat} that~$\dagger$, when restricted to the set of ASMs or the set of DPPs in Theorem~\ref{DRT},
provides a bijection  to the corresponding set in which $k_1$ and~$k_2$ are interchanged.  Therefore,
$Z^\ASM_n(x,y,z_1,z_2)$ and $Z^\DPP_n(x,y,z_1,z_2)$ are symmetric in $z_1$ and $z_2$, as already noted in~\eqref{ASMDPPint}.

The question of whether $\ast$ and $\dagger$ for ASMs can correspond to $\ast$ and $\dagger$ for DPPs under a bijection
between $\ASM(n)$ and $\DPP(n)$ will now be addressed.  Since $\ast$ and $\dagger$ are involutions, whose orbits are therefore of size~1 or~2,
this question reduces to a question of whether the numbers of $\ast$-invariant or $\dagger$-invariant objects are equal for $\ASM(n)$ and for $\DPP(n)$.

In the case of $\ast$, it follows from a result of de Gier, Pyatov and Zinn-Justin~\cite[Prop.~3, first equation]{DegPyaZin09},
initially conjectured by Mills, Robbins and Rumsey~\cite[Conj.~3S]{MilRobRum83},
that the sizes of the $\ast$-invariant sets
$\{A\in\ASM(n)\mid A^\ast=A,\;\nu(A)=p,\;\mu(A)=m,\;\rho_1(A)=k_1,\;\rho_2(A)=k_2\}$
and
$\{D\in\DPP(n)\mid D^\ast=D,\;\nu(D)=p,\;\mu(D)=m,\;\rho_1(D)=k_1,\;\rho_2(D)=k_2\}$
are equal for any~$n$,~$p$,~$m$,~$k_1$ and~$k_2$.
(It can be seen, using~\eqref{aststat} and the structure of $\ast$-invariant ASMs and DPPs, that each of these sets
is empty unless~$n$ is odd, $2p+m=n(n-1)/2$, $k_1=k_2=(n-1)/2$, and $m/2-(n-1)/4$ is a nonnegative integer.)
It now follows from this result and Theorem~\ref{DRT} that there
exist bijections between $\ASM(n)$ and $\DPP(n)$ with the properties that,
for each $A\in\ASM(n)$ and $D\in\DPP(n)$ which correspond under the bijection,~$A^\ast$ and~$D^\ast$ also correspond under the bijection,
$\nu(A)=\nu(D)$, $\mu(A)=\mu(D)$, $\rho_1(A)=\rho_1(D)$ and $\rho_2(A)=\rho_2(D)$.

In the case of $\dagger$, the last equation of~\eqref{dagstat} implies that if an ASM or DPP $X$ is $\dagger$-invariant then
$\rho_1(X)=\rho_2(X)$.  Furthermore, it can be seen using the definitions of $\dagger$
that, for $n\ge4$, there exist non-$\dagger$-invariant $n\times n$ ASMs~$A$ which satisfy $\rho_1(A)=\rho_2(A)$,
but that for all $n$, $D\in\DPP(n)$ is~$\dagger$-invariant if and only if $\rho_1(D)=\rho_2(D)$.  Therefore,
for $n\ge4$ there does not exist a bijection between $\ASM(n)$ and $\DPP(n)$ with the property that,
for each $A\in\ASM(n)$ and $D\in\DPP(n)$ which correspond under the bijection,~$A^\dagger$ and~$D^\dagger$ also correspond under the bijection.
However, this is not surprising, since $\dagger$ for an ASM involves its bulk structure, whereas~$\dagger$ for a DPP
involves only its top two rows.

\subsection{Further multiply-refined enumeration results}\label{doubrefsec}
In this section, other enumerative results for ASMs, and in some cases DPPs, involving one or more boundary statistics, but
mostly no bulk statistics,
are discussed briefly.

An ASM contains four natural boundary statistics corresponding to the positions of the~$1$'s in its first and last row and column.
More specifically, for an ASM $A$, these statistics can be taken as $\rho_1(A)$, \ldots, $\rho_4(A)$, where $\rho_1(A)$ and $\rho_2(A)$
are given by~\eqref{rho1A}--\eqref{rho2A}, and
\begin{align}\notag\rho_3(A)&=\text{number of 0's above the 1 in the first column of }A,\\[1mm]
\label{rho34A}\rho_4(A)&=\text{number of 0's below the 1 in the last column of }A.\end{align}

It can easily be seen, using the natural symmetry operations on ASMs, that the joint distribution on $\ASM(n)$ of
the bulk statistics~$\nu$ and~$\mu$ (as given by~\eqref{nuA}--\eqref{muA}) together with a single boundary statistic $\rho_i$
is independent of $i$, for $i=1,\ldots,4$.
It follows from the result~\eqref{MRR}, together with~\eqref{Zsing}
and~\eqref{ASMDPPint}, that these distributions on $\ASM(n)$ also equal the
joint distribution on $\DPP(n)$ of~$\nu$,~$\mu$ and~$\rho_j$ (as given by~\eqref{nuD}--\eqref{rho2D}), for $j=1,2$.
A determinant formula for these six identical distributions is provided by~\eqref{singdet}.

Now consider singly-refined enumeration without bulk statistics, and define
\begin{equation}\mathcal{A}_{n,k}=|\{A\in\ASM(n)\mid\rho_i(A)=k\}|=|\{D\in\DPP(n)\mid\rho_j(D)=k\}|,\end{equation}
for any $1\le i\le4$ and $1\le j\le2$.  Then
\begin{equation}\label{ASMDPPref}\mathcal{A}_{n,k}=
\textstyle\frac{(n+k-1)!\,(2n-k-2)!}{(2n-2)!\:k!\:(n-k-1)!}\;\prod_{l=0}^{n-2}\!\frac{(3l+1)!}{(n+l-1)!},\end{equation}
where this formula for the DPP case follows from results of Mills, Robbins and Rumsey~\cite[Sec.~5]{MilRobRum82}
(see Bressoud~\cite[Conj.~9 and Sec.~5.3]{Bre99} for further details), while the formula for the ASM case was first proved by
Zeilberger~\cite{Zei96b} (see Bressoud~\cite[Sec.~7.3]{Bre99} for further details),
following conjectures of Mills, Robbins and Rumsey~\cite[Conjs.~1 \& 2]{MilRobRum82},~\cite[Conjs.~1 \& 2]{MilRobRum83}.
Alternative proofs of the ASM case have been given by
Colomo and Pronko~\cite[Sec.~5.3]{ColPro05a},~\cite[Sec.~4.2]{ColPro06}, Fischer~\cite{Fis07}, and Stroganov~\cite[Sec.~4]{Str06}.
Each of these proofs of the ASM case of~\eqref{ASMDPPref} uses the six-vertex model with DWBC and particular forms of
the Izergin--Korepin formula, except for the proof of Fischer~\cite{Fis07}, which makes essential use of
a certain operator formula obtained by Fischer~\cite{Fis06,Fis10}.  Now define $\mathcal{A}_n=|\ASM(n)|=|\DPP(n)|=\sum_{k=0}^{n-1}\mathcal{A}_{n,k}$.
Then, as discussed in Section~\ref{intro}, and as can be obtained from~\eqref{ASMDPPref},
\begin{equation}\label{ASM}\textstyle\mathcal{A}_n=\prod_{i=0}^{n-1}\frac{(3i+1)!}{(n+i)!}.\end{equation}

Proceeding to doubly-refined enumeration,
it can easily be seen that the joint distributions on $\ASM(n)$ of the four statistics~$\nu$,~$\mu$,~$\rho_1$ and~$\rho_2$,
and the four statistics~$\nu$,~$\mu$,~$\rho_3$ and~$\rho_4$ (as given by~\eqref{nuA}--\eqref{rho2A} and~\eqref{rho34A})
are equal (and that each is symmetric in its two boundary statistics),
and Theorem~\ref{DRT} states
that these distributions are also equal to the joint distribution on $\DPP(n)$ of~$\nu$,~$\mu$,~$\rho_1$ and~$\rho_2$
(as given by~\eqref{nuD}--\eqref{rho2D}).
A determinant formula for these three identical distributions is provided by~\eqref{doubdet}.

Now consider doubly-refined enumeration without bulk statistics, and define
\begin{multline}\mathcal{A}_{n,i,j}=|\{A\in\ASM(n)\mid\rho_1(A)=i,\;\rho_2(A)=j\}|=|\{A\in\ASM(n)\mid\rho_3(A)=i,\;\rho_4(A)=j\}|\\
=|\{D\in\DPP(n)\mid\rho_1(D)=i,\;\rho_2(D)=j\}|.\end{multline}
It follows from~\eqref{propeq1} that
\begin{multline}\label{ASMdoubrefrec}
(\mathcal{A}_{n,i-1,j}-\mathcal{A}_{n,i,j-1})\,\mathcal{A}_{n-1}=\mathcal{A}_{n,i-1}\,\mathcal{A}_{n-1,j-1}-
\mathcal{A}_{n,i}\,\mathcal{A}_{n-1,j-1}-\\
\mathcal{A}_{n-1,i-1}\,\mathcal{A}_{n,j-1}+\mathcal{A}_{n-1,i-1}\,\mathcal{A}_{n,j},\end{multline}
which can be solved to give
\begin{multline}\label{ASMdoubref}
\mathcal{A}_{n,i,j}=\frac{1}{\mathcal{A}_{n-1}}\!\sum_{k=0}^{\min(i,n-j-1)}\bigl(\mathcal{A}_{n,i-k}\,\mathcal{A}_{n-1,j+k}-
\mathcal{A}_{n,i-k-1}\,\mathcal{A}_{n-1,j+k}-\\[-3mm]
\mathcal{A}_{n-1,i-k-1}\,\mathcal{A}_{n,j+k+1}+
\mathcal{A}_{n-1,i-k-1}\,\mathcal{A}_{n,j+k}\bigr).\end{multline}
Thus,~\eqref{ASMdoubref}, together with~\eqref{ASMDPPref}--\eqref{ASM}, provides an explicit formula for $\mathcal{A}_{n,i,j}$.
The ASM case of~\eqref{ASMdoubrefrec} was first obtained by Stroganov~\cite[Eq.~(34)]{Str06},
and the ASM case of~\eqref{ASMdoubref}, in a form also involving the ASM bulk statistics $\nu$ and~$\mu$, was obtained
by Colomo and Pronko~\cite[Eq.~(5.29)]{ColPro05b},~\cite[Eq.~(3.31)]{ColPro06}.

An expression for $Z^\ASM_n(1,1,z_1,z_2)=Z^\DPP_n(1,1,z_1,z_2)=\sum_{i,j=0}^{n-1}\mathcal{A}_{n,i,j}\,z_1^i\,z_2^j$
in terms of a particular Schur function can be obtained by considering~\eqref{ZCZASM} with the parameter~$q$ set to~$e^{\pm i\pi/3}$
or~$e^{\pm 2i\pi/3}$,
and using a formula of Okada~\cite[Thm.~2.4(1), 2nd equation]{Oka06} for the partition
function of the six-vertex model with DWBC at such $q$.
(See Di Francesco and Zinn-Justin~\cite[Eqs.~(2.2) \& (2.4)]{DifZin05}
or Biane, Cantini and Sportiello~\cite[Eq.~(1.5)]{BiaCanSpo12} for further details.)  It was shown by
Biane, Cantini and Sportiello~\cite[Thm.~1]{BiaCanSpo12}, using this expression together with a
determinantal identity for certain Schur functions,
which is itself obtained using a general identity involving minors of a matrix,
that the numbers $\mathcal{A}_{n,i,j}$ satisfy
$\det_{0\le i,j\le n-1}(\mathcal{A}_{n,i,j})=(-1)^{n(n+1)/2+1}\,(\mathcal{A}_{n-1})^{n-3}$.

The doubly-refined enumeration of ASMs without bulk statistics, and without reference to DPPs, has also been considered
in the context of totally symmetric self-complementary plane partitions (TSSCPPs).
In particular, it was conjectured by
Mills, Robbins and Rumsey~\cite[Conj.~3]{MilRobRum86} that the joint distributions of
certain pairs of statistics on TSSCPPs in a $2n\times2n\times2n$ box are equal to the joint
distribution of $\rho_1$ and $\rho_2$ (or $\rho_3$ and $\rho_4$) on $\ASM(n)$
(see also Robbins~\cite[p.~16]{Rob91}), and this has been proved
by Fonseca and Zinn-Justin~\cite{FonZin08}.
Pfaffian and constant-term expressions for the associated doubly-refined TSSCPP generating function
have been obtained by Ishikawa~\cite[Thm.~1.4, Cor.~7.3~\& Cor.~8.2]{Ish06a}, and,
using the result of~\cite{FonZin08} and Theorem~\ref{DRT},
these expressions, as well as certain integral expressions of~\cite[Eqs.~(4.9) \&~(4.14)]{FonZin08},
also apply to $Z^\ASM_n(1,1,z_1,z_2)=Z^\DPP_n(1,1,z_1,z_2)$.

Furthermore, the numbers $\mathcal{A}_{n,i,j}$ have been considered in the
context of a doubly-refined Razumov--Stroganov conjecture by Di Francesco~\cite[Sec.~4]{Dif04b}.

Proceeding now to other types of doubly-refined ASM enumeration,
for $i\ne j$ and $1\le i,j\le 4$, define $\mathcal{Z}^{i,j}_n(x,y,z_1,z_2)=
\sum_{A\in\ASM(n)}x^{\nu(A)}y^{\mu(A)}\,z_1^{\rho_i(A)}\,z_2^{\rho_j(A)}$,
where it can be checked that this is symmetric in $i$ and $j$.
Thus, for the already-considered case, which involves opposite boundaries of each ASM,
$Z^\ASM_n(x,y,z_1,z_2)=\mathcal{Z}^{1,2}_n(x,y,z_1,z_2)=\mathcal{Z}^{3,4}_n(x,y,z_1,z_2)$.
It can also be checked, using symmetry operations on $\ASM(n)$ and~\eqref{aststat}--\eqref{dagstat}, that for the remaining cases,
which involve adjacent boundaries of each ASM,
$\mathcal{Z}^{1,3}_n(x,y,z_1,z_2)=\mathcal{Z}^{2,4}_n(x,y,z_1,z_2)=
x^{n(n-1)/2}\,(z_1z_2)^{n-1}\,\mathcal{Z}^{1,4}_n(\tfrac{1}{x},\tfrac{y}{x},\tfrac{1}{z_1},\tfrac{1}{z_2})
=x^{n(n-1)/2}\,(z_1z_2)^{n-1}\,\mathcal{Z}^{2,3}_n(\tfrac{1}{x},\tfrac{y}{x},\tfrac{1}{z_1},\tfrac{1}{z_2})$.
No DPP statistic is currently known whose enumerative behavior together with
the DPP boundary statistic~$\rho_1$ (or $\rho_2$) matches that of the statistics associated with two adjacent boundaries of an ASM.
However, a remarkably simple relation between the generating functions for opposite-boundary and adjacent-boundary doubly-refined
enumeration of ASMs,
for the case without bulk statistics (i.e., $x=y=1$), has been obtained by Stroganov~\cite[p.~61]{Str06}.

Proceeding to quadruply-refined ASM enumeration, by following
an approach similar to
that of Section~\ref{ASMpropsect}, but using the form~\eqref{DesJac}, instead of~\eqref{Bazin}, of the Desnanot--Jacobi identity,
an identity which recursively determines
$\sum_{A\in\ASM(n)}x^{\nu(A)}y^{\mu(A)}\,z_1^{\rho_1(A)}\,z_2^{\rho_2(A)}\,z_3^{\rho_3(A)}\,z_4^{\rho_4(A)}$
has been obtained by Behrend~\cite[Thm.~1]{Beh12}.  An expression for this quadruply-refined ASM generating function
in the case $x=y=1$ has also been obtained by Ayyer and Romik~\cite[Thm.~2]{AyyRom12}.

Finally, multiply-refined ASM enumerations involving the configurations
of several rows or columns closest to ASM boundaries, but not involving bulk statistics, have been studied by
Fischer~\cite{Fis11,Fis12}, Fischer and Romik~\cite{FisRom09}, and Karklinsky and Romik~\cite{KarRom10}.
For example, in~\cite{Fis11,KarRom10}
a relation between opposite-boundary doubly-refined ASM enumeration and
doubly-refined ASM enumeration involving the configurations of
the first and second (or last and second-last) rows (or columns) of an ASM is obtained,
and a simple formula for the latter is derived, while in~\cite[Thm.~1]{Fis12} a relation between
two types of triply-refined ASM enumeration, one of which involves any three of the
statistics $\rho_1,\ldots,\rho_4$, is obtained. No multiply-refined enumerations of DPPs are currently known
which match the multiply-refined enumerations of ASMs in these cases.

\let\oldurl\url
\makeatletter
\renewcommand*\url{%
        \begingroup
        \let\do\@makeother
        \dospecials
        \catcode`{1
        \catcode`}2
        \catcode`\ 10
        \url@aux
}
\newcommand*\url@aux[1]{%
        \setbox0\hbox{\oldurl{#1}}%
        \ifdim\wd0>\linewidth
                \strut
                \\
                \vbox{%
                        \hsize=\linewidth
                        \kern-\lineskip
                        \raggedright
                        \strut\oldurl{#1}%
                }%
        \else
                \hskip0pt plus\linewidth
                \penalty0
                \box0
        \fi
        \endgroup
}
\makeatother
\gdef\MRshorten#1 #2MRend{#1}%
\gdef\MRfirsttwo#1#2{\if#1M%
MR\else MR#1#2\fi}
\def\MRfix#1{\MRshorten\MRfirsttwo#1 MRend}
\renewcommand\MR[1]{\relax\ifhmode\unskip\spacefactor3000 \space\fi
  \MRhref{\MRfix{#1}}{{\tiny \MRfix{#1}}}}
\renewcommand{\MRhref}[2]{%
 \href{http://www.ams.org/mathscinet-getitem?mr=#1}{#2}}
\bibliography{Bibliography}
\bibliographystyle{amsplainhyper}
\end{document}